\documentclass[12pt]{article}
\usepackage{amssymb,amsmath}
\usepackage{mathrsfs}
\usepackage{graphicx}

\newtheorem{theorem}{Theorem}
\newtheorem{lemm}{Lemma}
\newtheorem{corollary}{Corollary}

\begin{document}

\begin{center}
\large{\bf FLUCTUATIONS OF PROPAGATION FRONT\\
IN CATALYTIC BRANCHING WALK}
\end{center}
\vskip0,5cm
\begin{center}
Ekaterina Vl. Bulinskaya\footnote{ \emph{Email address:} {\tt
bulinskaya@yandex.ru}}$^,$\footnote{The work is supported by Russian Science Foundation under grant 17-11-01173 and is fulfilled
at Novosibirsk State University. The author is Associate Professor of the Lomonosov Moscow State University.}
\end{center}
\vskip1cm

\begin{abstract}
We consider a supercritical catalytic branching random walk (CBRW) on a multidimensional lattice $\mathbb{Z}^d$, $d\in\mathbb{N}$. The main subject of study is the behavior of particles cloud in space and time. For CBRW on an integer line, Carmona and Hu (2014) examined the asymptotical behavior of the maximal coordinate $M_n$ of the particles at time $n$. They proved that $M_n/n\to\mu$ almost surely (on a set of local non-degeneracy of CBRW), as $n\to\infty$, where $\mu>0$ is a certain constant. Under additional assumption of a single catalyst in CBRW they also investigated the fluctuations of $M_n$ with respect to $\mu n$, as $n\to\infty$.
Bulinskaya (2018) extended the strong limit theorem by Carmona and Hu having estimated the rate of the population propagation for the front of a multidimensional CBRW. Now our aim is to analyze fluctuations of the propagation front in CBRW on $\mathbb{Z}^d$. We not only solve the problem in a multidimensional setting but also, treating the case of an arbitrary finite number of catalysts for $d = 1$, generalize the result by Carmona and Hu  with the help of other probabilistic-analytic methods.


\vskip0,5cm 2010 {\it AMS classification}: 60J80, 60F05.

\end{abstract}

\section{Introduction}\label{s:intro_fluctuations}

Almost 50 years ago the study of population spread was initiated for a standard (or space-homogeneous) branching random walk (BRW) on $\mathbb{R}^d$, $d\in\mathbb{N}$, see, e.g., \cite{Biggins_78}. The problem proved difficult and has not been fully resolved. Its solution depends on quite a number of factors: characteristics of the random walk such as ``heaviness'' of the distribution tail of the walk jump,  characteristics of the branching, i.e. regime of the offspring reproduction (supercritical, critical or subcritical), the space dimension (in one-dimensional case the population is located between the minimum and the maximum of the particles coordinates, but in a multidimensional case the coordinates are ordered partially and one has to consider a random cloud of particles), the very setting of the problem (besides revealing the limit of properly normalized particles positions on the ``front'' one can study evolution of the particles positions oscillations around the limit, i.e. the front fluctuations) and etc. The contributions to this domain are due to many researchers, see, e.g., survey \cite{Shi_LNM_15}, recent papers \cite{Bhattacharya_etal_18}, \cite{Mallein_16} and others.

We are interested in a \emph{space-inhomogeneous}, namely \emph{catalytic}, BRW (CBRW) on an integer lattice $\mathbb{Z}^d$, $d\in\mathbb{N}$. Such a model was proposed back in the 90s of the 20th century, see, e.g., \cite{ABY_98}, \cite{HTV_12}, \cite{VTY} and references therein. Usually in these papers it is assumed that on the lattice $\mathbb{Z}^d$ there is a single catalyst (the source of branching, i.e. reproduction and death of particles), where particles may produce offspring and die, whereas outside it they perform random walk until the next hitting the catalyst. The main subjects of study were the total and local particles numbers at time $t$ and their behavior, as $t\to\infty$. In 2012 in the paper \cite{Molchanov_Yarovaya_12} for the first time there was set and solved the problem of the spread of particles population in supercritical CBRW, when the symmetric random walk has ``light'' tails and the particles reproduction is binary. There was also studied the rate of movement of the population propagation front, determined in terms of the moments boundedness of the local particles numbers, as $t\to\infty$. The study of the population spread in supercritical CBRW from another view-point of the almost sure convergence was initiated in paper \cite{Carmona_Hu_14} for $d=1$ and for ``light'' distribution tails of the random walk. This investigation was developed in the series of the author's works \cite{B_UMN_19}, \cite{B_SPA_18}--\cite{B_Innopolis_19} under assumptions of $d>1$ and ``heavy'' distribution tails of the walk jump.

In \cite{Carmona_Hu_14} it is established that the maximum $M_n$ at time $n$, i.e. the position of the right-most particle in CBRW on $\mathbb{Z}$, obeys the strong law of large numbers. In other words, $\lim_{n\to\infty}M_n/n=\mu>0$ a.s. (almost surely) on the set of local non-degeneracy of CBRW, where $\mu$ is a known constant. There is also found a non-trivial limit distribution for $M_n-\mu n$ as $n\to\infty$, i.e. there are studied the fluctuations of $M_n$ with respect to linear growth of $\mu n$. In the strong limit theorem the number of catalysts was assumed finite whereas in the weak limit theorem a \emph{single} catalyst was considered in view of technical difficulties. Our aim in the present work is to extend the results of paper \cite{Carmona_Hu_14} related to fluctuations to the case of $d>1$ and to the case of an \emph{arbitrary} finite catalysts number.

Both problems are solved. Namely, a method of measuring oscillations of the population front in CBRW on $\mathbb{Z}^d$ with respect to the scaled limit form of the front, at each point of this surface, is proposed and the corresponding limit distribution, as time grows to infinity, is found. As a consequence, for $d=1$ the result of the paper \cite{Carmona_Hu_14} is proved in extended setting for an arbitrary catalysts number in CBRW.

It is interesting to note the differences in the behavior of the maximum $M_n$ in a standard BRW on $\mathbb{Z}$ (see, e.g., \cite{Mallein_16}) and in CBRW considered in the present work. In both models the main term of the asymptotic behavior of $M_n$ is a linear function of time $n$. Besides the linear term in the representation of $M_n$, BRW has a negative logarithmic term and stochastically bounded fluctuations, whereas in the case of CBRW there are only stochastically bounded fluctuations. The same observation applies to a multidimensional model of BRW (see, e.g., \cite{Biggins_78}) versus CBRW on $\mathbb{Z}^d$, as show the results of this paper.

The methods we develop differ from those employed in the papers related to the study of BRW and from the methods of the works \cite{Carmona_Hu_14}, \cite{Molchanov_Yarovaya_12}. The investigation of asymptotic behavior of the maximum in one-dimensional case and the front in a multidimensional case is reduced to analysis of the solution to a non-linear integral equations system derived by us. Here we essentially bear on our previous results in papers \cite{B_TVP_14} and \cite{B_DAN_15}. Moreover, we use the apparatus of multidimensional renewal theory, theory of large deviations for sums of random variables, the Laplace transform, the Laplace method of studying the asymptotic behavior of Laplace integrals, auxiliary Bellman-Harris branching processes, hitting times under taboo and other probabilistic and analytic tools.

The structure of the work is as follows. In section~\ref{s:results_fluctuations} we give a description of the model and formulate the main results. Section~\ref{s:proofs_fluctuations} contains proofs. Primarily, the proofs are implemented under the assumption of a single catalyst located at the origin when $d=1$. Further on they are generalized to comprise an arbitrary catalysts number and the case of $d>1$, and that is the advantage of our approach.

\section{Description of CBRW and main results}\label{s:results_fluctuations}

All the random elements under consideration are supposed to be defined on the complete probability space $(\Omega,\mathcal{F},{\sf P})$, where $\Omega$ is a sample space consisting of possible outcomes $\omega\in\Omega$. Index ${\bf x}$ in expressions of the form ${\sf E}_{\bf x}$ and ${\sf P}_{\bf x}$ denotes the initial point of either CBRW or the random walk ${\bf S}$ depending on the context. Bold font of ${\bf x}$ emphasizes that ${\bf x}$ is a vector, whereas a symbol $x$ means that $x$ takes scalar values.

Recall the description of CBRW on $\mathbb{Z}^d$, $d\in\mathbb{N}$ (in our setting given in \cite{B_SPA_18}). At the initial time $t=0$ there is a single particle, which moves on $\mathbb{Z}^d$ according to a continuous-time Markov chain ${\bf S}=\{{\bf S}(t),t\geq0\}$, generated by the infinitesimal matrix ${Q=(q({\bf x},{\bf y}))_{{\bf x},{\bf y}\in\mathbb{Z}^d}}$. Assume that the underlying random walk (i.e. CBRW without branching) is homogeneous, the Markov chain ${\bf S}$ is irreducible, with conservative matrix $Q$, i.e. the matrix $Q$ has finite elements and
\begin{equation}\label{condition1}
q({\bf x},{\bf y})=q({\bf x}-{\bf y},{\bf 0}),\quad\sum\limits_{{\bf y}\in\mathbb{Z}^d}{q({\bf x},{\bf y})}=0,
\end{equation}
where $q({\bf x},{\bf y})\geq0$ for ${\bf x}\neq{\bf y}$ and $q:=-q({\bf x},{\bf x})\in(0,\infty)$, for all ${\bf x},{\bf y}\in\mathbb{Z}^d$. The requirement of the matrix $Q$ symmetry from the papers \cite{Molchanov_Yarovaya_12} and \cite{Platonova_Ryadovkin_17} is not imposed here.

When the particle hits a finite catalysts set $W=\{{\bf w}_1,\ldots,{\bf w}_N\}$, $W\subset\mathbb{Z}^d$, say, at point ${\bf w}_k$, it spends there random time, distributed according to an exponential law with parameter $\beta_k>0$, where $\beta_k=q/(1-\alpha_k)$, $k=1,\ldots,N$. Then the particle either splits or leaves the point ${\bf w}_k$ with the corresponding probabilities $\alpha_k$ and $1-\alpha_k$ ($0\leq\alpha_k<1$). If the particle splits at point ${\bf w}_k$, it produces a random non-negative integer number $\xi_{k}$ of offsprings, located at the same point ${\bf w}_k$, and instantly dies. Whenever the particle leaves ${\bf w}_k$, it jumps to point ${\bf y}\neq{\bf w}_k$ with probability $(1-\alpha_k)q({\bf w}_k,{\bf y})/q$ and resumes its motion, governed by the Markov chain ${\bf S}$, or, possibly, produces offspring, if the catalyst is located at point ${\bf y}$. It is stipulated that all the new particles behave as independent copies of the parent particle.

Denote by $f_k(s):={\sf E}{s^{\xi_k}}$, $s\in[0,1]$, the probability generating function of the random variable $\xi_k$, $k=1,\ldots,N$. We employ a standard assumption of finiteness of a derivative $f_k'(1)$, i.e. finiteness of the value $m_k:={\sf E}{\xi_k}$, for any $k=1,\ldots,N$. Moreover, in the present work we suppose that a well-known logarithmic moment condition is satisfied (see, e.g., \cite{Sev_71}, Ch.~2, Section~3), i.e.
\begin{equation}\label{condition2}
{\sf E}\,\xi_k\ln\left(\xi_k+1\right)<\infty,\quad k=1,\ldots,N.
\end{equation}

Likewise classical branching processes (see, e.g., \cite{Sev_71} and \cite{Vatutin_book_09}), in accordance with the works \cite{B_TVP_14} and \cite{B_DAN_15} every CBRW can be classified as supercritical, critical or subcritical depending on the relationship between characteristics $m_k$, $k=1,\ldots,N$, and certain probabilities of finiteness of hitting times under taboo (see, e.g., \cite{B_SPL_14}). In critical and subcritical regimes the population degenerates locally with probability~$1$. In supercritical regime only the mean total and local particles numbers grow exponentially fast, as time tends to infinity. The rate of the exponential growth denoted by $\nu$ is traditionally called the \emph{Malthusian parameter}, and in supercritical regime one has $\nu>0$. The precise definition of this parameter can be found in Section~3 of paper \cite{B_TVP_14}. In the sequel we consider \emph{supercritical} CBRW only, since in the framework of other regimes the problem of detecting the speed of the population spread turns out to be ill-posed.

Further on by virtue of condition (\ref{condition1}) and Theorem~1.2 in book \cite{Bremaud_99}, we will consider the version of the random walk ${\bf S}$ such that
\begin{equation}\label{S(t)_representation}
{\bf S}(t)={\bf x}+\sum_{i=1}^{\Pi(t)}{\bf Y}^i,\quad t\geq0,
\end{equation}
where ${\bf x}\in\mathbb{Z}^d$ is the initial point of the random walk, $\Pi=\{\Pi(t),t\geq0\}$ is a Poisson process with the constant intensity $q$ and ${\bf Y}^i$ is the value of the $i$-th jump of the random walk, $i\in\mathbb{N}$. The random vectors ${\bf Y}^i=\left(Y^i_1,\ldots,Y^i_d\right)$, $i\in\mathbb{N}$, are independent, identically distributed with ${\sf P}({\bf Y}^1={\bf y})=q({\bf 0},{\bf y})/q$, ${\bf y}\in\mathbb{Z}^d$, ${\bf y}\neq{\bf 0}$, and do not depend on $\Pi(t)$, for any $t\geq0$.

Assume that the tails of the random walk ${\bf S}$ are ``light'', i.e. for any ${\bf s}\in\mathbb{R}^d$ the function
\begin{eqnarray}\label{condition3}
H({\bf s})&:=&\sum_{{\bf x}\in\mathbb{Z}^d}e^{\langle {\bf s},{\bf x}\rangle}q({\bf 0},{\bf x})=\sum_{{\bf x}\in\mathbb{Z}^d}\left(e^{\langle {\bf s},{\bf x}\rangle}-1\right)q({\bf 0},{\bf x})\\
&=&q\left(\sum_{{\bf x}\in\mathbb{Z}^d,\,{\bf x}\neq{\bf 0}}e^{\langle {\bf s},{\bf x}\rangle}\frac{q({\bf 0},{\bf x})}{q}-1\right)=q\left({\sf E}e^{\langle{\bf s},{\bf Y}^1\rangle}-1\right)\nonumber
\end{eqnarray}
is finite, where the operation $\langle\cdot,\cdot\rangle$ denotes the inner product of vectors. It is easy to verify that the Hessian of function $H$ is positive-definite and, consequently, $H$ is a convex function. Set $\mathcal{R}:=\{{\bf r}\in\mathbb{R}^d:H({\bf r})=\nu\}$.

Let $Z(t)$ be a (random) set of particles existing in CBRW at time ${t\geq0}$. For a particle $z\in Z(t)$, denote by ${\bf X}^z(t)=\left(X^z_1(t),\ldots,X^z_d(t)\right)$ its position at time $t$. Consider the set
$$
\mathcal{I}:=\left\{\omega:\limsup_{t\to\infty}\{z\in Z(t):{\bf X}^z(t)\in W\}\neq\varnothing\right\}\in\mathcal{F}.
$$
To avoid operations with a continuum number of sets $\left\{A_t\right\}_{t\geq0}$ we just put $\limsup_{t\to\infty}A_t:=\cup_{m=1}^{\infty}\cap_{k=1}^{\infty}\cup_{n=k}^{\infty}A_{n/2^m}$, i.e. we deal with binary-rational values of parameter $t$ instead of its all non-negative values. For each $\omega\in\mathcal{I}$, there is an increasing to infinity sequence of \emph{binary-rational} values $t^{br}_l(\omega)$, $l\in\mathbb{N}$, such that at each time $t^{br}_l(\omega)$ there are particles at the catalysts set $W$. The event consisting of possible outcomes $\omega$ for which there exists a similar sequence of \emph{any} (not only binary-rational) values $t^{any}_l(\omega)$, $l\in\mathbb{N}$, has the same probability ${\sf P}(\mathcal{I})$, and we may call the set $\mathcal{I}$ {\it the event of infinite number of catalysts visits}. The behavior of CBRW on the complement $\overline{\mathcal{I}}$ of this set is almost surely trivial. Indeed, for values $t\geq T_0(\omega)$ large enough either CBRW degenerates or CBRW forms a system of several random walks (without branching), starting at time $T_0$ from ${\bf X}^z(\omega,T_0)$, $z\in Z(T_0)$, respectively. The supercritical regime of CBRW guarantees that ${\sf P}(\mathcal{I})>0$ (Theorem~4 in paper \cite{B_DAN_15}).

In case $d=1$ let $M_t:=\max\{X^z(t):z\in Z(t)\}$ be the maximum of CBRW at time $t$, i.e. the position of the right-most particle existing in CBRW at time $t$. Recall that according to \cite{Carmona_Hu_14}, for supercritical CBRW on $\mathbb{Z}$, the strong limit theorem $M_t/t\to\mu$ a.s. on set $\mathcal{I}$, as $t\to\infty$, holds under conditions (\ref{condition1}) and (\ref{condition3}), where $\mu:=\nu/r$ and $r>0$ is such that $H(r)=\nu$. A similar result for CBRW on $\mathbb{Z}^d$ is established in \cite{B_SPA_18} and states that, for each ${\bf x}\in\mathbb{Z}^d$, the following equalities are true
\begin{eqnarray}
\!\!\!\!\!\!\!\!& &{\sf P}_{\bf x}\left(\omega:\forall\varepsilon>0\;\exists T_1=T_1(\omega,\varepsilon)\;\mbox{such that}\;\forall t\geq T_1\;\mbox{and}\;\forall z\in Z(t),\;{\bf X}^z(t)/t\notin\mathcal{O}_{\varepsilon}\right)=1,\label{B_SPA_result_1}\\
\!\!\!\!\!\!\!\!& &{\sf P}_{\bf x}\left(\left.\omega:\!\forall\varepsilon\in(0,\nu)\exists T_2=T_2(\omega,\varepsilon)\;\mbox{such that}\;\forall t\geq T_2\;\exists z\in Z(t),\;{\bf X}^z(t)/t\notin\mathcal{Q}_{\varepsilon}\right|\mathcal{I}\right)\!=\!1,\label{B_SPA_result_2}
\end{eqnarray}
where the sets $\mathcal{O}_{\varepsilon}:=\{{\bf x}\in\mathbb{R}^d:\langle {\bf x},{\bf r}\rangle>\nu+\varepsilon\;\;\mbox{for at least one}\;\;{\bf r}\in\mathcal{R}\}$, $\varepsilon>0$, and
$\mathcal{Q}_{\varepsilon}:=\{{\bf x}\in\mathbb{R}^d:\langle{\bf x},{\bf r}\rangle<\nu-\varepsilon\;\mbox{for all}\;{\bf r}\in\mathcal{R}\}$, $\varepsilon\in(0,\nu)$. Put $\mathcal{O}:=\mathcal{O}_0$, $\mathcal{Q}:=\mathcal{Q}_0$ and $\mathcal{P}:=\partial\mathcal{Q}=\partial\mathcal{O}$, where $\partial\mathcal{S}$ denotes the boundary of set $\mathcal{S}\subset\mathbb{R}^d$. Thus, in the multidimensional case a counterpart of the limit $\mu$ is a surface $\mathcal{P}\subset\mathbb{R}^d$, called \emph{the limiting shape of the front} of the particles population in CBRW. Note that each set $\mathcal{Q}_{\varepsilon}$, $\mathcal{Q}$ or $\mathcal{P}\cup\mathcal{Q}$ is convex as an intersection of half-spaces (see, e.g., Theorem~2.1 in the monograph \cite{Rockafellar_73}). In paper \cite{B_SPA_18} there is also shown that the set $\mathcal{P}$ can be defined as
\begin{eqnarray}\label{P_def}
\mathcal{P}\!\!\!&=&\!\!\!\{{{\bf x}\in\mathbb{R}^d:{\langle{\bf x},{\bf r}\rangle\!=\!\nu}}\;\mbox{for a single value}\;{{\bf r}\in\mathcal{R}}\;\mbox{and}\;{\langle{\bf x},{\bf r}\rangle\!<\!\nu}\;\mbox{for other}\;{{\bf r}\!\in\!\mathcal{R}}\}\nonumber\\
\!\!\!& &\!\!\!\mbox{or as}\quad\mathcal{P}=\{{\bf z}({\bf r}):{\bf r}\in\mathcal{R}\},\quad\mbox{where}\quad{\bf z}({\bf r})=\nu\nabla H({\bf r})/\langle\nabla H({\bf r}),{\bf r}\rangle.
\end{eqnarray}

In paper \cite{Carmona_Hu_14} the fluctuations of the maximum $M_n$ were studied with respect to $\mu n$, as $n\to\infty$, in CBRW on $\mathbb{Z}$ with a single catalyst. In other words, there was found the rate of convergence in the strong limit theorem for $M_n$ and it turned out that $M_n-\mu n$ has a non-trivial limit distribution. Now we extend these results to the case of a multidimensional lattice and an arbitrary finite catalysts number. The notion of maximal position of a particle on a multidimensional lattice is not defined, and likewise it is not obvious how to consider fluctuations of the particles positions at time $t$ around the scaled limit surface.

We propose the following approach. In view of formula (\ref{P_def}), for each point ${\bf x}\in\mathcal{P}$, there exists a single value of the parameter ${\bf r}={\bf r}({\bf x})\in\mathcal{R}$ such that $\langle{\bf x},{\bf r}({\bf x})\rangle=\nu$, whereas, for other points ${\bf y}\in\mathcal{P}$, it is valid that $\langle{\bf y},{\bf r}({\bf x})\rangle<\nu$. It follows that ${\bf r}({\bf x})/|{\bf r}({\bf x})|$ is the normal vector to the surface $\mathcal{P}$ at point ${\bf x}$. We propose to measure the ``magnitude'' of fluctuations of the particles positions with respect to the level ${\bf x}t$ for the limit point ${\bf x}\in\mathcal{P}$ by projecting the coordinates of ${\bf X}^z(t)$ of a particle $z\in Z(t)$ onto the normal vector ${\bf r}({\bf x})/|{\bf r}({\bf x})|$ to the limiting surface at point ${\bf x}$. Let ${\bf v}\in\mathcal{P}$ and $z_{\bf v}\in Z(t)$ denote such a sequence of particles that ${\bf X}^{z_{\bf v}}(t)/t\to {\bf v}$ a.s. on set $\mathcal{I}$, as $t\to\infty$ (such a sequence always exists on set $\mathcal{I}$ in view of Theorem~2 in \cite{B_SPA_18}). Then for the particles sequence $z_{\bf y}\in Z(t)$ corresponding to ${\bf y}\neq{\bf x}$ one has $\langle{\bf X}^{z_{\bf y}}(t),{\bf r}({\bf x})\rangle/t\to\langle{\bf y},{\bf r}({\bf x})\rangle<\nu$, $t\to\infty$, a.s. on set $\mathcal{I}$. Moreover, $\langle{\bf X}^{z_{\bf x}}(t),{\bf r}({\bf x})\rangle/t\to\langle{\bf x},{\bf r}({\bf x})\rangle=\nu$, $t\to\infty$, a.s. on set $\mathcal{I}$. Therefore, the main contribution to the asymptotic behavior of $M_t({\bf r})-\nu t$ is due to the particles sequence $z_{\bf x}\in Z(t)$ rather than $z_{\bf y}\in Z(t)$, ${\bf y}\neq{\bf x}$, as $t\to\infty$, where $M_t({\bf r}):=\max\left\{\langle{\bf X}^z(t),{\bf r}\rangle:z\in Z(t)\right\}$.

Thus, our approach consisting in the study of the limit behavior of $M_t({\bf r})-\nu t$, as $t\to\infty$, and each ${\bf r}\in\mathcal{R}$, has two merits. Firstly, we get a way to measure relative oscillations of the front of the particles population at each point of the limiting shape of the front of the population propagation. Secondly, fluctuations of the front at other points of its limiting shape do not influence the random variable under consideration.

In statements of the main results there will arise the following function $\varphi(\lambda;{\bf x})$, $\lambda\geq0$, ${\bf x}\in\mathbb{Z}^d$, which has already appeared earlier in Theorem~4 of paper \cite{B_DAN_15} as the Laplace transform of the limit distribution of the normalized total and local particle numbers in CBRW on $\mathbb{Z}^d$. For ${\bf x}\in W$, this function is defined as a solution to the system of integral equations
\begin{eqnarray}\label{varphi(lambda,wj)=system_equations}
\!\!\!& &\varphi(\lambda;{\bf w}_j)=\alpha_j\int\nolimits_0^{\infty}{f_j(\varphi(\lambda e^{-\nu u};{\bf w}_j))\,dG_j(u)}\\
\!\!\!&+&\!\!(1-\alpha_j)\sum_{k=1}^N{\int\nolimits_0^{\infty}{\varphi(\lambda e^{-\nu u};{\bf w}_k)\,dG_{j,k}(u)}}
+(1-\alpha_j)\!\left(\!1\!-\!\sum_{k=1}^N{_{W_k}F_{{\bf w}_j,{\bf w}_k}(\infty)}\!\right),\; j=1,\ldots,N.\!\nonumber
\end{eqnarray}
For ${\bf x}\in\mathbb{Z}^d\setminus W$, the function $\varphi(\lambda;{\bf x})$, $\lambda\geq0$, admits the representation
\begin{equation}\label{varphi(lambda,wj)=system_equations+}
\varphi(\lambda;{\bf x})=\sum_{k=1}^N{\int\nolimits_0^{\infty}{\varphi(\lambda e^{-\nu u};{\bf w}_k)\,d {_{W_k}F_{{\bf x},{\bf w}_k}(u)}}}+1-\sum_{k=1}^N{_{W_k}F_{{\bf x},{\bf w}_k}(\infty)}.
\end{equation}
Here $G_j(t):=1-e^{-\beta_jt}$, $t\geq0$. The symbol $\ast$ denotes the convolution operation, whereas $G_{j,k}(t):=G_j\ast{_{W_k}\overline{F}_{{\bf w}_j,{\bf w}_k}(t)}$, $t\geq0$, $j,k=1,\ldots,N$. In its turn, the function ${_{W_k}\overline{F}_{{\bf w}_j,{\bf w}_k}(t)}$, $t\geq0$, is a cumulative distribution function (c.d.f.) of the first hitting time of point ${\bf w}_k$ by the random walk ${\bf S}$ after exit out of the starting point ${\bf w}_j$ under the taboo on states $W_k:=W\setminus\{{\bf w}_k\}$, whereas $_{W_k}F_{{\bf x},{\bf w}_k}(t)$, $t\geq0$, is a c.d.f. of the first hitting time of point ${\bf w}_k$ by the random walk ${\bf S}$ under the taboo on states $W_k$, when the starting point of the walk ${\bf S}$ is ${\bf x}$ (for hitting times under taboo, see, e.g., \cite{B_SPL_14}).

According to Lemma~7 of paper \cite{B_Arxiv_18}, under conditions (\ref{condition1}) and (\ref{condition2}) the system (\ref{varphi(lambda,wj)=system_equations}) has a unique solution $\varphi(\,\cdot\,;{\bf w}_j)$, $j=1,\ldots,N$, in the function class $\mathcal{C}_{\theta}$, for each $\theta=\left(\theta_1,\ldots,\theta_N\right)$, $\theta_i>0$, $i=1,\ldots,N$. Here the function classes are
$$
\mathcal{C}:=\left\{\left(\varphi(\,\cdot\,;{\bf w}_1),\ldots,\varphi(\,\cdot\,;{\bf w}_N)\right):\varphi(\,\cdot\,;{\bf w}_i)\;
\mbox{отображает}\;[0,\infty)\;\mbox{in}\;(0,1],\phantom{\frac{1}{2}}\right.
$$
$$
\left.\varphi(0;{\bf w}_i)=1\;\mbox{and}\;\lim_{\lambda\to0+}\frac{1-\varphi(\lambda;{\bf w}_i)}{\lambda}
>0,\;i=1,\ldots,N\right\},
$$
$$
\mathcal{C}_{\theta}:=\left\{\left(\varphi(\,\cdot\,;{\bf w}_1),\ldots,\varphi(\,\cdot\,;{\bf w}_N)\right)\in\mathcal{C}:
\lim_{\lambda\to0+}\frac{1-\varphi(\lambda;{\bf w}_i)}{\lambda}=\theta_i,\;i=1,\ldots,N\right\}\!,
$$
where $\theta=\left(\theta_1,\ldots,\theta_N\right)$, $\theta_i>0$, $i=1,\ldots,N$.

Note that, for $d=1$, the random variable $M_t$ takes integer values and the subtracted linear correction $\mu t$ takes real values. Hence, in corollary~\ref{C:d=1} below there arises a correction term $\{\mu t+y\}$, where, as usual, $\{s\}\in[0,1)$ is a fractional part of number $s\geq0$, $[s]$ is its integer part and $s=[s]+\{s\}$. Similarly, for $d>1$, the random variable $M_t({\bf r})$ has a lattice distribution, whenever all the pairwise ratios of coordinates of the vector ${\bf r}$ are rational numbers. Namely, if $r_i=\bar{r}_ir^{\ast}$, where $\bar{r}_i\in\mathbb{Z}$, and the greatest common divisor of all $\bar{r}_i$, $i=1,\ldots,d$, is $1$, then $M_t({\bf r})$ takes values of the form $r^{\ast}k$, $k\in\mathbb{Z}$ (a similar conclusion see, e.g., in \cite{Borovkov_book_09}, Lemma~13.3.1, p.~402). Conversely, if at least one relation of coordinates of the vector ${\bf r}$ occurs irrational, then the range of values of the random variable $M_t({\bf r})$ has ``concentration points'' and, consequently, the corresponding distribution is not a lattice one. This explains the necessity of introduction of a correction function $\chi(t;y):=r^{\ast}\{\nu t/r^{\ast}+y/r^{\ast}\}$, $t\geq0$, $y\in\mathbb{R}$, in the first case, whereas, for the sake of convenience, in the second case we set $\chi(t;y)=0$ for all $t\geq0$, $y\in\mathbb{R}$.

Now we are ready to formulate the main result of the paper.

\begin{theorem}\label{T:fluctuation_light}
Let conditions \emph{(\ref{condition1})}, \emph{(\ref{condition2})} and \emph{(\ref{condition3})} be satisfied for supercritical CBRW on $\mathbb{Z}^d$ with the Malthusian parameter $\nu$. Then, for each ${\bf x}\in\mathbb{Z}^d$, ${\bf r}\in\mathcal{R}$ and $y\in\mathbb{R}$, one has
\begin{equation}\label{main_result_fluctuation_light_d}
\lim_{t\to\infty}{\left({\sf P}_{\bf x}\left(M_t({\bf r})-\nu t\leq y\right)-\varphi\left(e^{-y+\chi(t;y)};{\bf x}\right)\right)}=0,
\end{equation}
where the function $\varphi(\lambda;{\bf x})$, $\lambda\geq0$, ${\bf x}\in\mathbb{Z}^d$, is a solution to the equations system \emph{(\ref{varphi(lambda,wj)=system_equations})} and \emph{(\ref{varphi(lambda,wj)=system_equations+})}. Moreover, the function $\varphi(\lambda;{\bf x})$ tends to the probability $1-{\sf P}_{\bf x}(\mathcal{I})$ of the local extinction of the population in CBRW, when $\lambda\to\infty$, for each fixed ${\bf x}\in\mathbb{Z}^d$.
\end{theorem}

It is important that the limit function $\varphi$ in Theorem~\ref{T:fluctuation_light} is defined uniquely, since a solution to the system (\ref{varphi(lambda,wj)=system_equations}) is searched for in the function class $\mathcal{C}_{\theta}$, where the vector $\theta=(\theta_1,\ldots,\theta_N)$ has coordinates $\theta_j>0$  equal to $\lim_{y\to+\infty}\lim_{t\to\infty}\left(e^{y-\chi(t;y)}{\sf P}_{{\bf w}_j}(M_t({\bf r})>\nu t+y)\right)$. The precise values of the latter expressions are found below in Lemma~\ref{L:limlim=theta_fluctuations} for a single catalyst (it is denoted by $c_{\ast}$) and its counterpart in case of several catalysts and a multidimensional lattice.

The statement of the theorem implying that the function $\varphi(\lambda;{\bf x})$ tends to the probability $1-{\sf P}_{\bf x}(\mathcal{I})$ of the local extinction of the population in CBRW, when $\lambda\to\infty$, for each fixed ${\bf x}\in\mathbb{Z}^d$, follows from Theorem~4 in paper \cite{B_DAN_15}. It means that we obtain a complete description of the fluctuations of the propagation front of the particles population in CBRW under its local non-degeneracy. Conversely, if the population degenerates locally (with probability $1-{\sf P}_{\bf x}(\mathcal{I})$), then, as noted above, either CBRW degenerates or it constitutes a system of several random walks (without branching) from some time moment. In both cases the study of ``the population front'' is out of the question.

As a consequence we write the result related to investigation of the maximum of CBRW on $\mathbb{Z}$, which was established in \cite{Carmona_Hu_14} under the assumption that there is a single catalyst located at the origin (we do not impose such restrictions).

\begin{corollary}\label{C:d=1}
If conditions \emph{(\ref{condition1})}, \emph{(\ref{condition2})} and \emph{(\ref{condition3})} are valid for supercritical CBRW on $\mathbb{Z}$ with the Malthusian parameter $\nu$, then, for each $x\in\mathbb{Z}$ and $y\in\mathbb{R}$, the following relation is true
\begin{equation}\label{main_result_fluctuation_light}
\lim_{t\to\infty}{\left({\sf P}_{x}\left(M_t-\mu t\leq y\right)-\varphi\left(e^{-ry+\{\mu t+y\}};x\right)\right)}=0.
\end{equation}
\end{corollary}

The proof of the main results differs essentially from the arguments of paper \cite{Carmona_Hu_14}, although in both works the renewal theory plays the key role. Whereas in \cite{Carmona_Hu_14} there are used estimates from above and from below for the probability under consideration, in our work we derive an equation for this probability and then find the asymptotic behavior of its solution. An advantage of our approach consists in that the proof, implemented initially for the case of a single catalyst and the lattice dimension $1$, is naturally extended to the case of many catalysts and an arbitrary lattice dimension.

\section{Proof of the main results}\label{s:proofs_fluctuations}

For the sake of exposition clarity, at first consider CBRW on $\mathbb{Z}$ with a single catalyst $w_1$ located, without loss of generality, at the origin, i.e. $W=\{w_1\}$ with $w_1=0$, and the starting point~$0$ as well.

Let $E(t;u):={\sf P}_0\left(\exists z\in Z(t):X^z(t)>u\right)={\sf P}_0\left(M_t>u\right)$, $t,u\geq0$. The following lemma proved in \cite{B_Arxiv_18} contains an integral equation for the probability $E(t;u)$.

\begin{lemm}\label{L:equation_multi}
Let condition \emph{(\ref{condition1})} be valid. Then the probability $E(t;u)$, $t,u\geq0$, satisfies the non-linear integral equation of the convolution type
\begin{equation}\label{E(t;u)_equation}
E(t;u)=\alpha_1\int\nolimits_0^t{\left(1-f_1\left(1-E(t-s;u)\right)\right)\,dG_1(s)}
+(1-\alpha_1)\int\nolimits_0^t{E(t-s;u)\,dG_{1,1}(s)}+I\left(t;u\right),
\end{equation}
where
\begin{eqnarray}\label{I(t;u)_definition}
I(t;u)&:=&{\sf P}_0\left(S(t)>u\right)-\int\nolimits_0^t{{\sf P}_0\left(S(t-s)>u\right)\,d F_{0,0}(s)}\\
&-&\alpha_1\int\nolimits_0^t{{\sf P}_0\left(S(t-s)>u\right)\,d\left(G_1(s)-G_1\ast F_{0,0}(s)\right)}.\nonumber
\end{eqnarray}
Here the function $F_{0,0}(t)$, $t\geq0$, is a c.d.f. of the first hitting time of point $0$ by the random walk $S$, when the starting point of $S$ is $0$. Similarly, the function $\overline{F}_{0,0}(t)$, $t\geq0$, is a c.d.f. of the first hitting time of point $0$ by the random walk $S$ after exit out of the starting point $0$ and $G_{1,1}(t):=G_1\ast\overline{F}_{0,0}(t)$, $t\geq0$.
\end{lemm}

In the following lemma there are auxiliary results related to probabilities of large deviations of the random walk under consideration (without branching).

\begin{lemm}\label{L:P(S(t)>x)_from_discrete_to_continuous}
Let conditions \emph{(\ref{condition1})} and \emph{(\ref{condition3})} be satisfied for the random walk $S$. Then, for all ${s,t,x\geq0}$, the following inequality holds true
\begin{equation}\label{P(S(t)>x)_estimate1}
{\sf P}_0\left(S(t)\geq x\right)\leq e^{-sx+tH(s)}.
\end{equation}
Moreover, for all $t\geq0$, $\theta:=x/t\geq q{\sf E}Y^1$, it is valid that
\begin{equation}\label{P(S(t)>x)_estimate2}
{\sf P}_0\left(S(t)\geq x\right)\leq e^{-t\Lambda(\theta)},
\end{equation}
where the function $\Lambda\left(\vartheta\right):=\sup_{s\in\mathbb{R}}\left(\vartheta s-\ln{\sf E}_0e^{sS(1)}\right)=\sup_{s\in\mathbb{R}}\left(\vartheta s-H(s)\right)$, $\vartheta\in\mathbb{R}$, is the deviation function of the random variable $S(1)$ \emph{(}see, e.g., \emph{\cite{Borovkov_Borovkov_08}}, Ch.\emph{~6}, Section\emph{~1)}.
If additionally $\theta=x/t>q{\sf E}Y^1$, $x\in\mathbb{Z}$, $t>0$, then for all such $x,t\to\infty$, uniformly in $\theta\in[q{\sf E}Y^1+\varepsilon_1,\Theta_1]$ for each $\varepsilon_1>0$ and $\Theta_1>q{\sf E}Y^1+\varepsilon_1$, one has
\begin{equation}\label{P(S(t)>x)_asymptotics}
{\sf P}_0\left(S(t)\geq x\right)\sim\frac{e^{-t\Lambda(\theta)}}{\left(1-e^{-\lambda(\theta)}\right)\sqrt{2\pi t D(\theta)}},
\end{equation}
where $\lambda\left(\vartheta\right):=\Lambda'\left(\vartheta\right)>0$, when $\vartheta>{\sf E}S(1)=q{\sf E}Y^1$, and $D\left(\vartheta\right):=\left.H''(s)\right|_{s=\lambda(\vartheta)}$, $\vartheta\in\mathbb{R}$, is a variance of some random variable.
\end{lemm}

{\sc Proof. }With the help of representation (\ref{S(t)_representation}) and the exponential Chebyshev's inequality, applied to a discrete-time random walk (see, e.g., inequality (1.1.19) in Theorem~1.1.1 in \cite{Borovkov_book_13}, Ch.~1, Section~1), in the case of continuous time one has
$$
{\sf P}_0\left(S(t)\geq x\right)\!=\!\sum_{j=0}^{\infty}{\sf P}\left(\Pi(t)\!=\!j\right){\sf P}_0\left(\sum\limits_{i=1}^j Y^i\geq x\right)\leq\sum_{j=0}^{\infty}e^{-qt}\frac{(qt)^j}{j!}e^{-s x}\left({\sf E}e^{sY^1}\right)^j=e^{-sx+tH(s)},
$$
for all $s,t,x\geq0$. Thus, statement (\ref{P(S(t)>x)_estimate1}) is proved. Inequality (\ref{P(S(t)>x)_estimate2}) follows from the proven inequality (\ref{P(S(t)>x)_estimate1}) and the definition of the function $\Lambda$.

Taking into account representation (\ref{S(t)_representation}) and the apparatus of characteristic functions, it is not difficult to verify that the random walk $S$ is a random process with independent increments. Moreover, similarly to the proof of inequality (\ref{P(S(t)>x)_estimate1}), we derive that ${\sf E}_0e^{sS(t)}=e^{tH(s)}$. Set $t_k:=[2^{k}t]/2^{k}$. For each fixed $k\in\mathbb{Z}_+$, we have $t_k\to\infty$ if and only if $t\to\infty$. Additionally, $t_k\leq t$ and $t-t_k\to0$, as $k\to\infty$. Employing asymptotic formula (6.1.17) from Corollary~6.1.7 in book \cite{Borovkov_Borovkov_08}, Ch.~6, Section~1, valid for a discrete-time random walk, in the framework of the lemma conditions we obtain
\begin{equation}\label{P(S(t)>x)_asymptotics_binary-rational}
{\sf P}_0\!\left(S\!\left(t_k\right)\!\geq\!x\right)\!=\!{\sf P}_0\!\left(\!\sum_{i=1}^{[2^k t]}\!\left(S\!\left(i/2^k\right)\!-\!S\!\left((i-1)/2^k\right)\right)\!\geq\! x\!\right)\!\!\sim\!\frac{e^{-[2^k t]\Lambda_{k}(x/[2^k t])}}
{\!\left(\!1\!-\!e^{-\lambda_{k}(x/[2^k t])}\!\right)\!\!\sqrt{2\pi [2^k t]D_{k}(x/[2^k t])}},
\end{equation}
for each fixed $k\in\mathbb{Z}_+$ and $x,t\to\infty$. Asymptotic relation (\ref{P(S(t)>x)_asymptotics_binary-rational}) holds uniformly in $x/[2^k t]\in[{\sf E}_0S\left(1/2^k\right)+\varepsilon_2,\Theta_2]$, where $\varepsilon_2$ is any positive number and the value $\Theta_2>{\sf E}_0S\left(1/2^k\right)$ is chosen arbitrarily as well. In formula (\ref{P(S(t)>x)_asymptotics_binary-rational}) the function $\Lambda_{k}(\vartheta):=\sup_{s\in\mathbb{R}}\left(s\vartheta-\ln{\sf E}_0e^{sS(1/2^k)}\right)$, $\vartheta\in\mathbb{R}$, is the deviation function of the random variable $S(1/2^k)$, $\lambda_k(\vartheta):=\Lambda_k'(\vartheta)$ and $D_k(\vartheta):=\left.\left(\ln{\sf E}_0e^{sS\left(1/2^k\right)}\right)''\right|_{s=\lambda_k(\vartheta)}$. It is easy to see that $\Lambda_{k}(\vartheta)=\sup_{s\in\mathbb{R}}\left(\vartheta s-H(s)/2^k\right)=\Lambda\left(2^k\vartheta\right)/2^k$, $\lambda_k(\vartheta)=\lambda(2^k\vartheta)$ and $D_k(\vartheta)=\left.H''(s)\right|_{s=\lambda\left(2^k\vartheta\right)}/2^k$, for all $\vartheta\in\mathbb{R}$ and each $k\in\mathbb{Z}_+$. Therefore, relation (\ref{P(S(t)>x)_asymptotics_binary-rational}) can be rewritten as follows
\begin{equation}\label{P(S(t)>x)_asymptotics_t_k}
{\sf P}\left(S(t_k)\geq x\right)\sim\frac{e^{-t_k\Lambda(\theta_k)}}{\left(1-e^{-\lambda(\theta_k)}\right)\sqrt{2\pi t_k D(\theta_k)}},\quad t,x\to\infty,
\end{equation}
for each $k\in\mathbb{Z}_+$ uniformly in $\theta_k:=x/t_k\in[q{\sf E}Y^1+\varepsilon_1,\Theta_1]$, whenever we consider $\varepsilon_2=\varepsilon_1/2^k$ and $\Theta_2=\Theta_1/2^k$. However, the uniform convergence in $\theta_k$ implies the uniform convergence in $k\in\mathbb{Z}_+$ as well. Thus, formula (\ref{P(S(t)>x)_asymptotics})
is established for $t=n/2^k$ and $n\to\infty$ uniformly in $k\in\mathbb{Z}_+$.

To complete the proof of Lemma~\ref{L:P(S(t)>x)_from_discrete_to_continuous} we have to pass from the case of binary-rational values of $t$ to arbitrary real ones. Recall that $\theta=x/t$ and $\theta_k=x/t_k$. Clearly, $t\geq t_k$ and $\theta\leq\theta_k$. Denote also $g(\hat{t},\hat{\theta}):=\hat{t}\Lambda(\hat{\theta})+\ln\hat{t}/2$ and $1/h(\hat{\theta}):=\left(1-e^{-\lambda(\hat{\theta})}\right)\sqrt{2\pi D(\hat{\theta})}$, $\hat{t}>0$, $\hat{\theta}>q{\sf E}Y^1$. Then, for all values of $t$ large enough and all $\theta,\theta_k\in[q{\sf E}Y^1+\varepsilon_1,\Theta_1]$, one has
\begin{eqnarray}\label{P(S(t)>x)_binary-rational-P(S(t)>x)}
\!\!\!\!\!\!\!\!& &\!\!\!\!\!\left|\frac{{\sf P}_0\left(S(t)\!\geq\!x\right)}{h(\theta)e^{-g(t,\theta)}}-\frac{{\sf P}_0\left(S\left(t_k\right)\!\geq\! x\right)}{h(\theta_k)e^{-g(t_k,\theta_k)}}\right|\!\leq\!\frac{\left|{\sf P}_0\left(S(t)\geq x\right)-{\sf P}_0\left(S\left(t_k\right)\geq x\right)\right|}{h(\theta)e^{-g(t,\theta)}}\!+\!\frac{{\sf P}_0\left(S\left(t_k\right)\geq x\right)}{h(\theta_k)e^{-g(t_k,\theta_k)}}\\
\!\!\!\!\!\!\!\!&\times&\!\!\!\!\!\frac{\left|h(\theta)e^{-g(t,\theta)}\!-\!h(\theta_k)e^{-g(t_k,\theta_k)}\right|}
{h(\theta)e^{-g(t,\theta)}}\!\leq\!\frac{q\sqrt{t}(t\!-\!t_k)}{e^{-t\Lambda(\Theta_1)}}
\max(h(\vartheta))^{\!-\!1}\!+\!C_1(t\!-\!t_k)\!+\!C_2e^{g(t,\theta)-g(t_k,\theta)}(t\!-\!t_k),\nonumber
\end{eqnarray}
where $\max\left(h(\vartheta)\right)^{-1}$ is taken over values $\vartheta$ from the interval $[q{\sf E}Y^1+\varepsilon_1,\Theta_1]$, whereas the constants $C_i=C_i(q{\sf E}Y^1+\varepsilon_1,\Theta_1)$, $i=1,2$, do not depend on $\theta$ and $\theta_k$ by virtue of relation (\ref{P(S(t)>x)_asymptotics_t_k}). Indeed, the latter inequality holds true, since condition (\ref{condition1}) implies that
$$
\left|{\sf P}_0\left(S(t)\geq x\right)-{\sf P}_0\left(S\left(t_k\right)\geq x\right)\right|\leq{\sf P}_0\left(S(t-t_k)\neq0\right)\leq1-e^{-q\left(t-t_k\right)}\leq q\left(t-t_k\right)
$$
and, moreover, the following relations are valid
$$
\left|h(\theta)e^{-g(t,\theta)}-h(\theta_k)e^{-g(t_k,\theta_k)}\right|
=\left|\left(h(\theta)e^{-g(t,\theta)}-h(\theta_k)e^{-g(t,\theta)}\right)
+\left(h(\theta_k)e^{-g(t,\theta)}-h(\theta_k)e^{-g(t,\theta_k)}\right)\right.
$$
$$
\left.+\left(h(\theta_k)e^{-g(t,\theta_k)}-h(\theta_k)e^{-g(t_k,\theta_k)}\right)\right|\leq e^{-g(t,\theta)}\left(\theta_k-\theta\right)\max h'(\vartheta)\!+\!
h(\theta_k)e^{-g(t,\theta)}\!\left(g(t,\theta_k)\!-\!g(t,\theta)\right)\!
$$
$$
+h(\theta_k)e^{-g(t_k,\theta_k)}\left(g(t,\theta_k)-g(t_k,\theta_k)\right)\leq e^{-g(t,\theta)}\theta\frac{t-t_k}{t_k}\max h'(\vartheta)+
h(\theta_k)e^{-g(t,\theta)}t\left(\theta_k-\theta\right)\Lambda'(\theta_1)
$$
$$
+h(\theta_k)e^{-g(t_k,\theta_k)}\left((t-t_k)\Lambda(\theta_k)+\ln(t/t_k)/2\right)\leq e^{-g(t,\theta)}\Theta_1\frac{t-t_k}{t-1}\max h'(\vartheta)+h(\theta_k)e^{-g(t,\theta)}\frac{t}{t-1}$$
$$
\times(t-t_k)\Theta_1\lambda(\Theta_1)+h(\theta_k)e^{-g(t_k,\theta)}(t-t_k)\Lambda(\Theta_1)+h(\theta_k)e^{-g(t_k,\theta)}\frac{t-t_k}{t-1}.
$$
Here we use inequalities $1-e^{-u}\leq u$, $\ln u\leq u-1$, $u\geq0$, $g(t_k,\theta_k)\geq g(t_k,\theta)$, and the mean value theorem.

Finally, consider $k=k(t)$ and set, e.g., $k(t)=[t]^2$. Then $t-t_k\leq 2^{-k(t)}\to0$, as $t\to\infty$, and $e^{g(t,\theta)-g(t_k,\theta)}\to1$, as $t\to\infty$, uniformly in $\theta\in[q{\sf E}Y^1+\varepsilon_1,\Theta_1]$. Therefore, it follows from relation (\ref{P(S(t)>x)_binary-rational-P(S(t)>x)}) that
$$
\left|\frac{{\sf P}_0\left(S(t)\geq x\right)}{h(\theta)e^{-g(t,\theta)}}-\frac{{\sf P}_0\left(S\left(t_k\right)\geq x\right)}{h(\theta_k)e^{-g(t_k,\theta_k)}}\right|\to0,\quad x,t\to\infty,
$$
uniformly in $\theta,\theta_k\in[q{\sf E}Y^1+\varepsilon_1,\Theta_1]$, which in combination with relation (\ref{P(S(t)>x)_asymptotics_t_k}) implies the required formula (\ref{P(S(t)>x)_asymptotics}).

Lemma~\ref{L:P(S(t)>x)_from_discrete_to_continuous} is proved completely. $\square$

The definition of the supercritical regime of CBRW (see \cite{B_TVP_14}) entails two following formulae ${\alpha_1m_1+(1-\alpha_1)F_{0,0}(\infty)>1}$ and $\alpha_1m_1G_1^{\ast}(\nu)+(1-\alpha_1)\,G_1^{\ast}(\nu)\overline{F}^{\,\ast}_{0,0}(\nu)=1$. In terms of the function $G(t):=\alpha_1m_1G_1(t)+(1-\alpha_1)\,G_1\ast\overline{F}_{0,0}(t)$, $t\geq0$, it means that $G^{\ast}(\nu)=1$. Here $J^{\ast}(\lambda)$, $\lambda\geq0$, is the Laplace transform of a c.d.f. $J(t)$, $t\geq0$, with the support on non-negative semi-axes, i.e. $J^{\ast}(\lambda):=\int\nolimits_{0-}^{\infty}{e^{-\lambda t}\,dJ(t)}$.

\begin{lemm}\label{L:fluctuations_main_step}
Let conditions \emph{(\ref{condition1})} and \emph{(\ref{condition3})} be satisfied. Then, for each fixed ${y\in\mathbb{R}}$, one has
$$
\lim_{t\to\infty}{e^{-r\{\mu t+y\}}\int_0^t{I(t-u;\mu t+y)\,d\sum_{k=0}^{\infty}{G^{\ast k}(u)}}}=c_{\ast}e^{-ry},
$$
where the constant
$$
c_{\ast}:=\frac{e^{-r}\left(1-F_{0,0}^{\ast}(\nu)\right)\left(1-\alpha_1G_1^{\ast}(\nu)\right)}
{\sqrt{2}(1-e^{-r})H'(r)\int\nolimits_0^{\infty}{se^{-\nu s}\,dG(s)}}.
$$
\end{lemm}

{\sc Proof.} The function $\Lambda(\vartheta)$, $\vartheta\in\mathbb{R}$, is the Fenchel-Legendre transform of the function $H(s)$, $s\in\mathbb{R}$, (see, e.g., \cite{Borovkov_Borovkov_08}, Ch.~6, Section~1). In particular, it means that $\Lambda'(\vartheta)=\lambda(\vartheta)=\lambda$ if and only if $H'(\lambda)=\vartheta$, for all $\vartheta$, $\lambda\in\mathbb{R}$, and $\Lambda(\vartheta)=\lambda(\vartheta)\vartheta-H(\lambda(\vartheta))$. Recall that $H(r)=\nu$ and denote $H'(r)=\theta_0$. Consequently, $\Lambda'(\theta_0)=\lambda(\theta_0)=r$ and $\Lambda(\theta_0)=r\theta_0-\nu$. In view of the convexity of the function $H$ combined with equalities $H(0)=0$ and $H(r)=\nu$, one has $H'(0)<\nu/r=\mu$, where $H'(0)={\sf E}_0S(1)=q{\sf E}Y^1$. Since the function $H'(s)$, $s\in\mathbb{R}$, is increasing one, $\theta_0=H'(r)>H'(0)$. The function $\Lambda(\vartheta)$, $\vartheta\in\mathbb{R}$, attains its minimal value, equal to $0$, at point $\vartheta={\sf E}_0S(1)=H'(0)$, and, therefore, $\Lambda(\theta_0)>0$. In other words, $r\theta_0-\nu>0=r\mu-\nu$, i.e. $\theta_0>\mu$.

Show, for each $y\in\mathbb{R}$ and any $\varepsilon\in(0,\theta_0-\mu)$, validity of the estimate
\begin{equation}\label{int=o}
\int\nolimits_{\Gamma(t)}{I(t-u;\mu t+y)\,d\sum_{k=0}^{\infty}{G^{\ast k}(u)}}=o(1),\quad t\to\infty,
\end{equation}
where $\Gamma(t):=[0,t]\setminus[t-t\mu(\theta_0-\varepsilon)^{-1},t-t\mu(\theta_0+\varepsilon)^{-1}]\subset\mathbb{R}$, $t\geq 0$. The proofs of Lemmas~1 and 2 in paper~\cite{B_Arxiv_18} imply relations $0\leq I(t;u)\leq{\sf P}_0\left(S(t)>u\right)$, for all $t,u\geq0$, under condition (\ref{condition1}). Using inequality (\ref{P(S(t)>x)_estimate1}) we get
\begin{equation}\label{int_o_leq}
\int\nolimits_{\Gamma(t)}{I(t-u;\mu t+y)\,d\sum_{k=0}^{\infty}{G^{\ast k}(u)}}
\leq\int\nolimits_{\Gamma(t)}{{\sf P}_0(S(t-u)>\mu t+y)\,d\sum_{k=0}^{\infty}{G^{\ast k}(u)}}=:J_1(t)+J_2(t),
\end{equation}
for all $t\geq0$. The latter integral is the sum of two integrals $J_1(t)$ and $J_2(t)$, related to the integration intervals $\Gamma(t)\cap[0,t-v(t)]$ and $\Gamma(t)\cap(t-v(t),t]$, respectively, where $v(t)$, $t\geq0$, is an increasing to infinity function, $v(t)=o(t)$, $t\to\infty$.

To estimate the integral $J_1(t)$ apply inequality (\ref{P(S(t)>x)_estimate2}) to the integrand for values $t$ large enough, which results in
\begin{eqnarray}\label{J_1(t)<=}
{\sf P}_0\left(S(t-u)>\mu t+y\right)\!&\!\leq\!&\!\exp\left\{-(t-u)\Lambda\left(\frac{\mu t+y}{t-u}\right)\right\}=\exp\left\{-(t-u)\left(r\cdot\frac{\mu t+y}{t-u}-\nu\right)\right\}\nonumber\\
\!&\!\times\!&\!\exp\left\{-(t-u)V(t,u)\right\}=\exp\left\{-ry-\nu u-(t-u)V(t,u)\right\},
\end{eqnarray}
where the function $V(t,u):=\Lambda\left((\mu t+y)/(t-u)\right)-r\left(\mu t+y\right)/(t-u)+\nu$, $t>u$, $u\geq0$, takes non-negative values only by virtue of the definition of the function $\Lambda$. Moreover, ${V(t,u)=0}$ if and only if $\Lambda'\left((\mu t+y)/(t-u)\right)=\lambda\left((\mu t+y)/(t-u)\right)=r$. The latter equality takes place for $(\mu t+y)/(t-u)=\theta_0$, i.e. for $u=t-(\mu t+y)/\theta_0$. Thus, if $|\theta_0y/(\mu t+y)|<\varepsilon$ (this is valid for all large enough values of $t$), then $u\notin\Gamma(t)$. Hence, the function $V(t,u)$ does not vanish for $u\in\Gamma(t)$ and all large enough $t$. Since the function $\Lambda(\vartheta)-r\vartheta+\nu$, $\vartheta\in\mathbb{R}$, is convex, it attains its minimal value at point $\vartheta=r$. It follows that, for any $\varepsilon$ and large enough $t$, there exists $\delta>0$ such that $V(t,u)\geq\delta$, for all $u\in\Gamma(t)$. On account of formula~(\ref{J_1(t)<=}) and Theorem~25 in book \cite{Vatutin_book_09}, p.~30, we come to the estimate
\begin{eqnarray}\label{J_1(t)=o}
\!\!\!\!\!\!\!\!\!\!\!&\!\!&\!J_1(t)\leq e^{-ry}\!\!\!\!\!\int\limits_{\Gamma(t)\cap[0,t-v(t)]}\!\!\!\!\!{e^{-\nu u-(t-u)\delta}\,d\sum_{k=0}^{\infty}{G^{\ast k}(u)}}\leq e^{-ry-(t-v(t))\nu-v(t)\delta}\nonumber\\
\!\!\!\!\!\!\!\!\!\!\!&\!\times\!&\!\int\limits_0^{t-v(t)}{e^{(t-v(t)-u)(\nu-\delta)}\,d\sum_{k=0}^{\infty}{G^{\ast k}(u)}}
\sim e^{-ry-v(t)\delta}\left(\int\limits_0^{\infty}{se^{-\nu s}\,dG(s)}\right)^{-1}\int\limits_0^{\infty}{e^{-\delta s}\,ds}=o(1),
\end{eqnarray}
as $t\to\infty$, whenever $v(t)\to\infty$.

To estimate the integral $J_2(t)$ apply inequality (\ref{P(S(t)>x)_estimate1}), when $s=r+\varepsilon_0$ for some $\varepsilon_0>0$, and also Theorem~25 in book \cite{Vatutin_book_09}, p.~30. As a result we have
\begin{eqnarray}\label{J_2(t)_o}
\!\!\!\!\!J_2(t)&\leq&\int_{t-v(t)}^t{e^{-(r+\varepsilon_0)(\mu t+y)+(t-u)H(r+\varepsilon_0)}\,d\sum_{k=0}^{\infty}{G^{\ast k}(u)}}\leq e^{-(r+\varepsilon_0)y-\nu t-t\mu\varepsilon_0+v(t)H(r+\varepsilon_0)}\nonumber\\
\!\!\!\!\!&\times&\int_{0}^t{d\sum_{k=0}^{\infty}{G^{\ast k}(u)}}\leq e^{-(r+\varepsilon_0)y-t\mu\varepsilon_0+v(t)H(r+\varepsilon_0)}\left(\nu\int\nolimits_0^{\infty}{se^{-\nu s}dG(s)}\right)^{-1}=o(1),
\end{eqnarray}
as $t\to\infty$, when $v(t)=o(t)$. Combination of the established formulae (\ref{int_o_leq}), (\ref{J_1(t)=o}) and (\ref{J_2(t)_o}) implies the desired relation (\ref{int=o}).

To complete the proof of Lemma~\ref{L:fluctuations_main_step} we have to show that
\begin{equation}\label{int_neq_o}
\lim_{t\to\infty}{e^{-r\{\mu t+y\}}\int\nolimits_{t-t\mu(\theta_0-\varepsilon)^{-1}}^{t-t\mu(\theta_0+\varepsilon)^{-1}}{I(t-u;\mu t+y)\,d\sum_{k=0}^{\infty}{G^{\ast k}(u)}}}=c_{\ast}e^{-ry}.
\end{equation}
At first we study the asymptotic behavior, as $t\to\infty$, of the following integral
$$
\int\limits_{t-t\mu/(\theta_0-\varepsilon)}^{t-t\mu/(\theta_0+\varepsilon)}\!\!\!\!{{\sf P}_0\left(S(t-u)>\mu t+y\right)\!\!\,d\sum_{k=0}^{\infty}{G^{\ast k}(u)}}=-\!\!\!\!\int\limits_{\mu/(\theta_0+\varepsilon)}^{\mu/(\theta_0-\varepsilon)}\!\!\!{{\sf P}_0\left(S(tv)>\mu t+y\right)\,d\sum_{k=0}^{\infty}{G^{\ast k}(t(1-v))}},
$$
which we denote by $J_3(t)$ and where we have made the change of variables $v=1-u/t$, i.e. $u=t(1-v)$. Note that ${\sf P}_0\left(S(tv)>\mu t+y\right)={\sf P}_0\left(S(tv)\geq\mu t+y+1-\{\mu t+y\}\right)$. Therefore, employing formula (\ref{P(S(t)>x)_asymptotics}), we obtain
$$
J_3(t)\sim-\!\!\!\int\limits_{\mu/(\theta_0+\varepsilon)}^{\mu/(\theta_0-\varepsilon)}\!\!\!{\frac{e^{-tv\Lambda\left(\mu/v+\rho/v\right)}}
{\left(1-e^{-\lambda\left(\mu/v+\rho/v\right)}\right)\sqrt{2\pi tv D\left(\mu/v+\rho/v\right)}}\,d\sum_{k=0}^{\infty}{G^{\ast k}(t(1-v))}},\quad t\to\infty,
$$
where $\rho=\rho(t):=(y+1-\{\mu t+y\})/t$. Now analysis of the asymptotic behavior of the latter integral is reduced to the study of the following function
\begin{equation}\label{J_3(t)sim}
J_3(t)\sim\sqrt{t}\!\!\!\int\limits_{\mu/(\theta_0+\varepsilon)}^{\mu/(\theta_0-\varepsilon)}\!\!\!
{\frac{C_3e^{-t\left(v\Lambda\left(\mu/v+\rho/v\right)-\nu(1-v)\right)}}
{\left(1-e^{-\lambda\left(\mu/v+\rho/v\right)}\right)\sqrt{2\pi v D(\mu/v+\rho/v)}}\,dv},\quad t\to\infty.
\end{equation}
Indeed, taking into account the asymptotic behavior of the renewal density function (see, e.g., \cite{Cox_62}, p.~55) one has
\begin{equation}\label{renewal_function_sim}
d\sum_{k=0}^{\infty}{G^{\ast k}(u)}=\sum_{k=0}^{\infty}{g^{\ast k}(u)}\,du,\quad e^{-\nu t}\sum_{k=0}^{\infty}{g^{\ast k}(t)}=\sum_{k=0}^{\infty}{\tilde{g}^{\ast k}(t)},\quad \sum_{k=0}^{\infty}{\tilde{g}^{\ast k}(t)}\to C_3,\quad t\to\infty,
\end{equation}
where the constant $(C_3)^{-1}:=\int\nolimits_0^{\infty}{se^{-\nu s}dG(s)}=\int\nolimits_0^{\infty}{se^{-\nu s}g(s)\,ds}$, $g(u):=G\,'(u)$, $\tilde{g}(u):=e^{-\nu u}g(u)$, $u\geq0$, and $g^{\ast k}$ stands for the convolution of the $k$-th order of the density $g$.

The integral in relation (\ref{J_3(t)sim}) is a Laplace integral (see, e.g., \cite{Fedoryuk_87}, p.~96) containing an additional parameter $\rho$. Therefore, Theorem~2.1 in book \cite{Fedoryuk_87}, Ch.~2, Section~2, entails the following asymptotic relation
$$
J_3(t)\sim\frac{C_3e^{-t\rho r}}{\left(1-e^{-r}\right)\theta_0\sqrt{2D(\theta_0)\Lambda''(\theta_0)}},\quad t\to\infty.
$$
Since $\Lambda''(\vartheta)=\lambda'(\vartheta)$ and $H'(\lambda(\vartheta))=\vartheta$, i.e. $H''(\lambda(\vartheta))\lambda'(\vartheta)=1$, one has $\Lambda''(\theta_0)=1/H''(r)$, where $\theta_0=H'(r)$ and $\lambda(\theta_0)=r$. As a result we come to relation
$$
J_3(t)\sim\frac{e^{-r\left(y+1-\{\mu t+y\}\right)}}{\left(1-e^{-r}\right)H'(r)\sqrt{2}\int\nolimits_0^{\infty}{se^{-\nu s}\,dG(s)}},\quad t\to\infty.
$$

With the help of (\ref{renewal_function_sim}) we establish that
\begin{equation}\label{G_1_ast_sum_Gk_sim}
\int\nolimits_0^t{\sum_{k=0}^{\infty}{g^{\ast k}(t-u)}\,dF_{0,0}(u)}\sim C_3F_{0,0}^{\ast}(\nu)e^{\nu t},\quad
\int\nolimits_0^t{\sum_{k=0}^{\infty}{g^{\ast k}(t-u)}\,dG_1(u)}\sim C_3G_1^{\ast}(\nu)e^{\nu t},
\end{equation}
$$
\int\nolimits_0^t{\sum_{k=0}^{\infty}{g^{\ast k}(t)}\left(G_1\ast F_{0,0}(u)\right)}\sim C_3G_1^{\ast}(\nu)F_{0,0}^{\ast}(\nu)e^{\nu t},\quad t\to\infty,
$$
where $G_1^{\ast}(\nu)=\beta_1/\left(\nu+\beta_1\right)$. Hence, taking into account definition (\ref{I(t;u)_definition}) and operating similarly to deriving the asymptotic behavior  of the function $J_3(t)$, as $t\to\infty$, we get relation (\ref{int_neq_o}). Lemma~\ref{L:fluctuations_main_step} is proved completely. $\square$

In Lemma~\ref{L:fluctuations_main_step}, for each fixed $y\in\mathbb{R}$ and $t\to\infty$, we find the asymptotic behavior of the integral $\int_0^t{I(t-u;\mu t+y)\,d\sum_{k=0}^{\infty}{G^{\ast k}(u)}}$. The following lemma yields an estimate of this integral from above, for \emph{all} $y,t\geq0$.

\begin{lemm}\label{L:int<=}
Let conditions \emph{(\ref{condition1})} and \emph{(\ref{condition3})} be satisfied. Then, for all $y,t\geq0$, the inequality
\begin{equation}\label{int<=}
\int_0^t{I(t-u;\mu t+y)\,d\sum_{k=0}^{\infty}{G^{\ast k}(u)}}\leq C_4e^{-ry}
\end{equation}
is valid, where $C_4$ is a positive constant.
\end{lemm}
{\sc Proof.} The proof of Lemma~\ref{L:int<=} mainly resembles the proof of the previous Lemma~\ref{L:fluctuations_main_step}. Firstly, once again apply the estimate $0\leq I(t;u)\leq{\sf P}_0\left(S(t)>u\right)$, $t,u\geq0$, following from the proof of Lemmas~1 and 2 in paper \cite{B_Arxiv_18} under condition (\ref{condition1}). Secondly, use the following representation
\begin{equation}\label{intP=J1+J2+J3}
\int_0^t{{\sf P}_0\left(S(t-u)>\mu t+y\right)\,d\sum_{k=0}^{\infty}{G^{\ast k}(u)}}=J_1(t;y)+J_2(t;y)+J_3(t;y),
\end{equation}
where the integrals $J_1(t;y)$, $J_2(t;y)$ and $J_3(t;y)$, $t,y\geq0$, differ by the integration areas only which are $\Upsilon_1(t;y):=\left[t-(\mu t+y)(1-\varepsilon_3)/\theta_0,t-v(t)\right]$, $(t-v(t),t]$ and $[0,t-(\mu t+y)(1-\varepsilon_3)/\theta_0)$, respectively. Here a fixed positive number $\varepsilon_3\in(0,1)$, and $v(t)$, $t\geq0$, is an increasing to infinity function such that $v(t)=o(t)$, $t\to\infty$, and $v(t)\leq(\mu t+y)(1-\varepsilon_3)/\theta_0$, for all $t\geq0$. If there exists $\varepsilon_4>0$ such that $(\mu t+y)/\theta_0>t(1+\varepsilon_4)$, then we set $\Upsilon_1(t;y):=[0,t-v(t)]$ and $J_3(t;y)=0$.

Similarly to the derivation of the estimate for $J_1(t)$ in Lemma~\ref{L:fluctuations_main_step} (see formulae (\ref{J_1(t)<=}) and (\ref{J_1(t)=o})) we come to inequalities
\begin{equation}\label{J_1(t;y)<=}
J_1(t;y)\leq e^{-ry}\!\!\!\int\limits_{\Upsilon_1(t;y)}\!\!\!{e^{-\nu u-(t-u)V(t,u)}\,d\sum_{k=0}^{\infty}G^{\ast k}(u)}\leq e^{-ry}\!\!\int\nolimits_0^{t-v(t)}\!\!\!\!\!{e^{-\nu u-(t-u)\delta_1}\,d\sum_{k=0}^{\infty}G^{\ast k}(u)}\leq C_5e^{-ry},
\end{equation}
valid for all $y\geq0$, some positive constant $C_5$ and any $t\geq T_3$. Recall that a non-negative function $V(t,u)$, $t>u$, $u\geq0$, vanishes at point $u=t-(\mu t+y)/\theta_0$ only, which does not belong to the integration area $\Upsilon_1(t;y)$. Consequently, on set $u\in\Upsilon_1(t;y)$ either estimate $V(t,u)\geq V\left(t,t-(\mu t+y)(1-\varepsilon_3)/\theta_0\right)=\Lambda\left(\theta_0/(1-\varepsilon_3)\right)-r\theta_0/(1-\varepsilon_3)+\nu\geq\delta_1$ or estimate $V(t,u)\geq V(t,0)>\Lambda(\theta_0(1+\varepsilon_4))-r\theta_0(1+\varepsilon_4)+\nu\geq\delta_1$ holds true for some value $\delta_1$.

Similarly to derivation of the estimate for $J_2(t)$ in relation (\ref{J_2(t)_o}) in Lemma~\ref{L:fluctuations_main_step}, we obtain
\begin{equation}\label{J_2(t;y)<=}
J_2(t;y)\leq C_6e^{-ry},
\end{equation}
for all $y\geq0$, some positive constant $C_6$ and all $t\geq T_4$.

Let us estimate $J_3(t;y)$ from above in a non-trivial case when $J_3(t;y)\neq 0$. Employing relation (\ref{P(S(t)>x)_asymptotics}) we come to inequalities
\begin{eqnarray}\label{J_3(t;y)<=}
\!& &\!J_3(t;y)\leq C_7\int_0^{t-(\mu t+y)(1-\varepsilon_3)/\theta_0}{\frac{e^{-(t-u)\Lambda\left((\mu t+y)/(t-u)\right)}}{\sqrt{t-u}}\,d\sum_{k=0}^{\infty}G^{\ast k}(u)}\leq\frac{C_7\sqrt{\theta_0}e^{-ry}}{\sqrt{(\mu t+y)(1-\varepsilon_3)}}\nonumber\\
\!&\times&\!\int_0^{t-(\mu t+y)(1-\varepsilon_3)/\theta_0}{e^{-\nu u-(t-u)\Lambda''(\mu)\left((\mu t+y)/(t-u)-\theta_0\right)^2}\,d\sum_{k=0}^{\infty}G^{\ast k}(u)}\leq\frac{C_8e^{-ry}}{\sqrt{\mu t+y}}\\
\!&\times&\!\int_0^{t-(\mu t+y)(1-\varepsilon_3)/\theta_0}{e^{-\Lambda''(\mu)\left(\mu t+y-\theta_0(t-u)\right)^2/\sqrt{t-u}}\,du}\leq\frac{C_8e^{-ry}\sqrt{t}}{\theta_0\sqrt{\mu t+y}}\int_{-\infty}^{+\infty}{e^{-\Lambda''(\mu)v^2}\,dv}\leq C_9e^{-ry},\nonumber
\end{eqnarray}
valid for some constants $C_7$, $C_8$, $C_9$ (depending on $\varepsilon_3$) and all $t\geq T_5$. Here we use relation (\ref{renewal_function_sim}), the Taylor expansion in the form $\Lambda\left((\mu t+y)/(t-u)\right)=r\theta_0-\nu+r\left((\mu t+y)/(t-u)-\theta_0\right)+\Lambda''(\tilde{\theta})\left((\mu t+y)/(t-u)-\theta_0\right)^2$ (the value $\tilde{\theta}$ belongs to the interval between $(\mu t+y)/(t-u)$ and $\theta_0$) and the variable change $v=\left(\mu t+y-\theta_0(t-u)\right)/\sqrt{t}$, i.e. $dv=\theta_0du/\sqrt{t}$.

Finally, it follows from inequality (\ref{P(S(t)>x)_estimate1}) when $s=r$ that
$$
\int_0^t{{\sf P}_0\left(S(t-u)>\mu t+y\right)\,d\sum_{k=0}^{\infty}{G^{\ast k}(u)}}\leq e^{-ry}\int\nolimits_0^{T_6}{e^{-\nu u}\,d\sum_{k=0}^{\infty}{G^{\ast k}(u)}}=C_{10}e^{-ry},
$$
for a constant $C_{10}>0$ and all $t\in[0,T_6]$, where $T_6:=\max\{T_3,T_4,T_5\}$. The latter inequality combined with the established relations (\ref{intP=J1+J2+J3}), (\ref{J_1(t;y)<=}), (\ref{J_2(t;y)<=}) and (\ref{J_3(t;y)<=}) completes the proof of Lemma~\ref{L:int<=}, where we may set $C_4:=\max\{C_5,C_6,C_9,C_{10}\}$. $\square$

Next derive an estimate for the probability ${\sf P}_0\left(M_t-\mu(t+\tilde{t}\,)\geq y\right)$ from above.

\begin{lemm}\label{L:lower_estimate_fluctuations}
If conditions \emph{(\ref{condition1})} and \emph{(\ref{condition3})} are satisfied, then, for all $\tilde{t},t\geq0$, $y\in\mathbb{R}$ and some positive constant $C$ the following inequality is valid
\begin{equation}\label{P_0(Mt_>_Lt+r_u)_upper_estimate_fluctuations}
{\sf P}_0\left(M_t-\mu (t+\tilde{t}\,)>y\right)\leq Ce^{-\nu\tilde{t}-ry}.
\end{equation}
\end{lemm}
{\sc Proof. }For any $u\geq0$, according to the mean value theorem, applied to function $f_1$, equation (\ref{E(t;u)_equation}) entails the inequality
$$
E(t;u)\leq\int\nolimits_0^t{E(t-s;u)\,d G(s)}+I\left(t;u\right).
$$
Iterating this inequality $k$ times we get
$$E(t;u)\leq\int\nolimits_0^t{E(t-s;u)\,d G^{\ast (k+1)}(s)}+\int\nolimits_0^t{I(t-s;u)\,d\sum_{j=0}^k{G^{\ast j}(s)}}.$$
For any fixed $t$, one has $G^{\ast k}(t)\to 0$, as $k\to\infty$. This is valid, e.g., in view of Lemma~22 in book \cite{Vatutin_book_09}. Hence, the term $\int\nolimits_0^t{{\sf P}_0\left(M_{t-s}>u\right)\,dG^{\ast(k+1)}(s)}$ is negligibly small for large $k$. Therefore, the latter inequality can be rewritten in the form
\begin{equation}\label{P_0(Ltgamma/Mtgamma<lambda)_inequality_fluctuations}
E(t;u)\leq\int\nolimits_0^t{I(t-s;u)\,d\sum_{j=0}^{\infty}{G^{\ast j}(s)}}.
\end{equation}
Letting $u=\mu(t+\tilde{t}\,)+y$ in this relation and employing Lemma~\ref{L:int<=}, we come to the statement of Lemma~\ref{L:lower_estimate_fluctuations} when $\mu\tilde{t}+y\geq0$ and $C=C_4$. Whenever $\mu\tilde{t}+y<0$, the required inequality (\ref{P_0(Mt_>_Lt+r_u)_upper_estimate_fluctuations}) remains valid by taking $C\geq1$. Lemma~\ref{L:lower_estimate_fluctuations} is proved completely. $\square$

For $t\geq0$ and $u\geq0$, set
$$
J\left(t;u\right):=m_1\int\nolimits_0^t{E(t-s;u)\,
dG_1(s)}-\int\nolimits_0^t{\left(1-f_1\left(1-E(t-s;u)\right)\right)\,dG_1(s)}.
$$

\begin{lemm}\label{L:I-1(t;lambdaLt)ast_sim_fluctuations}
Whenever conditions \emph{(\ref{condition1})}, \emph{(\ref{condition2})} and \emph{(\ref{condition3})} hold true, the following relation takes place
$$
\lim_{y\to+\infty}\lim_{t\to\infty}e^{ry}\int\nolimits_0^t{J\left(t-s;\mu t+y\right)\,d\sum_{j=0}^{\infty}{G^{\ast j}(s)}}=0.
$$
\end{lemm}
{\sc Proof. }The mean value theorem, applied to function $f_1$, and Lemma~\ref{L:lower_estimate_fluctuations}, employed when $Ce^{-ry}<1$, guarantee that
$$
0\leq e^{ry}\int\nolimits_0^t{J\left(t-s;\mu t+y\right)\,d\sum_{j=0}^{\infty}{G^{\ast j}(s)}}
$$
$$
\leq C\int\nolimits_0^t{\left(m_1-f_1'\left(1-C e^{-r y} e^{-\nu s}\right)\right)e^{-\nu s}\,d\left(G_1\ast\sum_{j=0}^{\infty}{G^{\ast j}}(s)\right)}.
$$
In view of relation~(\ref{G_1_ast_sum_Gk_sim}) one has
$$
e^{ry}\int\nolimits_0^t{J(t-s;\mu t+y)\,d\sum_{j=0}^{\infty}{G^{\ast j}(s)}}\leq C_{11}\int\nolimits_0^t{\left(m_1-f_1'\left(1-Ce^{-ry} e^{-\nu s}\right)\right)\,ds},
$$
for some positive constant $C_{11}$. Let us show that the latter integral converges, as $t\to\infty$, whenever ${\sf E}\,\xi_1\ln{\xi_1}<\infty$. Indeed,
$$
\int\nolimits_0^t{\left(m_1-f_1'\left(1-Ce^{-ry}e^{-\nu s}\right)\right)\,ds}=\int\nolimits_0^t{\left({\sf E}\xi_1-{\sf E}\left(\xi_1\left(1-Ce^{-ry}e^{-\nu s}\right)^{\xi_1-1}\right)\right)\,ds}
$$
$$=\nu^{-1}{\sf E}\left(\xi_1\int\nolimits_{1-Ce^{-ry}}^{1-Ce^{-ry}e^{-\nu t}}{\frac{1-u^{\xi_1-1}}{1-u}\,du}\right)=\nu^{-1}{\sf E}\left(\xi_1\int\nolimits_{1-Ce^{-ry}}^{1-Ce^{-ry}e^{-\nu t}}{\sum_{k=1}^{\xi_1-1}u^{k-1}}\,du\right)
$$
$$=\nu^{-1}{\sf E}\!\left(\xi_1\!\sum_{k=1}^{\xi_1-1}\!\frac{\left(1-Ce^{-ry}e^{-\nu t}\right)^k-\left(1-Ce^{-ry}\right)^k}{k}\right)\!\leq\nu^{-1}{\sf E}\!\left(\xi_1\!\sum_{k=1}^{\xi_1-1}\!\frac{1-\left(1-Ce^{-ry}\right)^k}{k}\!\right)
$$
$$
\leq\nu^{-1}{\sf E}\left(\xi_1\left(1-\left(1-Ce^{-ry}\right)^{\xi_1}\right)\sum_{k=1}^{\xi_1-1}\frac{1}{k}\right)\leq\nu^{-1}{\sf E}\left(\xi_1\ln{\xi_1}\left(1-\left(1-Ce^{-ry}\right)^{\xi_1}\right)\right).
$$
Here we made the variable change $u=1-Ce^{-ry}e^{-\nu s}$, $du=C\nu e^{-ry}e^{-\nu s}\,ds$, i.e. $ds=du/\left(\nu(1-u)\right)$. Applying the Lebesgue bounded convergence theorem we see that, as $y\to+\infty$, one has ${\sf E}\left(\xi_1\ln{\xi_1}\left(1-\left(1-Ce^{-ry}\right)^{\xi_1}\right)\right)\to0$. Thus, Lemma~\ref{L:I-1(t;lambdaLt)ast_sim_fluctuations} is proved completely. $\square$

Lemma~\ref{L:lower_estimate_fluctuations} entails that $\lim_{y\to+\infty}E\left(t;\mu\left(t+\tilde{t}\,\right)+y\right)=0$ and $e^{ry}E\left(t;\mu\left(t+\tilde{t}\right)+y\right)\leq C e^{-\nu\tilde{t}}$, $y\in\mathbb{R}$, $t,\tilde{t}\geq0$. The next result refines the latter assertion with the help of Lemma~\ref{L:I-1(t;lambdaLt)ast_sim_fluctuations}, as $t\to\infty$ and afterwards $y\to+\infty$.

\begin{lemm}\label{L:limlim=theta_fluctuations}
If conditions \emph{(\ref{condition1})}, \emph{(\ref{condition2})} and \emph{(\ref{condition3})} hold true, then the following relation is valid
$$
\lim_{y\to+\infty}\lim_{t\to\infty}e^{ry-r\{\mu t+y\}}E\left(t;\mu t+y\right)=c_{\ast}.
$$
\end{lemm}
{\sc Proof. }In view of equation (\ref{E(t;u)_equation}), for any $u\geq0$, one has
$$
E\left(t;u\right)=\int\nolimits_0^t{E\left(t-s;u\right)\,d G(s)}+I(t;u)-J(t;u).
$$
Iteration of this equation $k$ times leads to
$$
E\left(t;u\right)\!=\!\int\nolimits_0^t{E\left(t-s;u\right)\,dG^{\ast(k+1)}(s)}\!+\!\int\nolimits_0^t{I(t-s;u)\,d\sum_{j=0}^k{G^{\ast j}(s)}}\!-\!\int\nolimits_0^t{J(t-s;u)\,d\sum_{j=0}^k{G^{\ast j}(s)}}.
$$
Applying Lemma~22 in \cite{Vatutin_book_09} once again, for each fixed $t$, we get $G^{\ast k}(t)\to 0$, as $k\to\infty$. Consequently, the term $\int\nolimits_0^t{E\left(t-s;u\right)\,dG^{\ast(k+1)}(s)}$ is negligibly small for large values of $k$. Therefore, the latter equation can be rewritten as follows
$$
E\left(t;\mu t+y\right)=\int\nolimits_0^t{I\left(t-s;\mu t+y\right)\,d\sum_{j=0}^{\infty}{G^{\ast j}(s)}}
-\int\nolimits_0^t{J\left(t-s;\mu t+y\right)\,d\sum_{j=0}^{\infty}{G^{\ast j}(s)}}.
$$
By dividing the both parts of the obtained equality by $e^{-ry}$, tending $t$ to infinity and afterwards $y$ to $+\infty$, we deduce the assertion of Lemma~\ref{L:limlim=theta_fluctuations} on account of Lemmas~\ref{L:fluctuations_main_step} and \ref{L:I-1(t;lambdaLt)ast_sim_fluctuations}. $\square$

The statement of the following lemma coincides with that of Corollary~\ref{C:d=1}, when $N=1$, $w_1=0$ and the starting point of CBRW is $x=0$.

\begin{lemm}\label{L:limK=0_fluctuations}
Let conditions \emph{(\ref{condition1})}, \emph{(\ref{condition2})} and \emph{(\ref{condition3})} be true. Then, for each $y\in\mathbb{R}$, we have
$$
\lim_{t\to\infty}{e^{ry-r\{\mu t+y\}}\left(1-E\left(t;\mu t+y\right)-\varphi(e^{-ry+r\{\mu t+y\}};0)\right)}=0,
$$
where $\varphi(\cdot;0)\in\mathcal{C}_{c_{\ast}}$.
\end{lemm}
{\sc Proof. }Let $K(t;y):=e^{ry-r\{\mu t+y\}}\left(1-E\left(t;\mu t+y\right)-\varphi(e^{-ry+r\{\mu t+y\}};0)\right)$. Firstly, note that
\begin{equation}\label{limlimK(t,lambda)=0_fluctuations}
\limsup_{y\to+\infty}\limsup_{t\to\infty}|K(t;y)|=0.
\end{equation}
It is valid by virtue of Lemma~\ref{L:limlim=theta_fluctuations}, equation (\ref{varphi(lambda,wj)=system_equations}) and the triangle inequality, since
$$
|K(t;y)|\leq\left|e^{ry-r\{\mu t+y\}}E\left(t;\mu t+y\right)-c_{\ast}\right|+\left|e^{ry-r\{\mu t+y\}}\left(1-\varphi(e^{-ry+r\{\mu t+y\}};0)\right)-c_{\ast}\right|.
$$

To prove the lemma, it is sufficient to verify that $K(y):=\lim_{T\to\infty}K_T(y)=0$, where $K_T(y):=\sup_{t\geq T}|K(t;y)|$. Equations (\ref{varphi(lambda,wj)=system_equations}) (when $N=1$ and $w_1=x=0$) and (\ref{E(t;u)_equation}) imply that
\begin{equation}\label{K(t;lambda)leqI2+I3_fluctuations}
e^{ry-r\{\mu t+y\}}\left(E\left(t;\mu t+y\right)-1+\varphi(e^{-ry+r\{\mu t+y\}};0)\right)=I_{11}(t,T;y)+I_{12}(t,T;y)+I_2(t;y),
\end{equation}
where for $T<t$ we set the function $I_{11}(t,T;y)$ to be 
\begin{eqnarray*}
\!\!\!& &\!\!\alpha_1e^{ry-r\{\mu t+y\}}\int\nolimits_0^{t-T}{\!\!\!\left(f_1(\varphi(e^{-ry+r\{\mu t+y\}}e^{-\nu s};0))-f_1\left(1-E\left(t-s;\mu t+y\right)\right)\right)\,dG_1(s)}\\
\!\!\!&+&\!\!(1-\alpha_1)e^{ry-r\{\mu t+y\}}\int\nolimits_0^{t-T}{\!\!\!\left(E\left(t-s;\mu t+y\right)-1+\varphi(e^{-ry+r\{\mu t+y\}}e^{-\nu s};0)\right)\,d\left(G_1\ast\overline{F}_{0,0}(s)\right)},
\end{eqnarray*}
whereas the integral $I_{12}(t,T;y)$ differs from $I_{11}(t,T;y)$ by the integration area only, i.e. there arises $\int_{t-T}^t$ instead of $\int_0^{t-T}$. Finally,
\begin{eqnarray*}
I_2(t;y)\!\!&:=&\!\!e^{ry-r\{\mu t+y\}}I\left(t;\mu t+y\right)\!-\!\alpha_1e^{ry-r\{\mu t+y\}}\int\nolimits_t^{\infty}\!\!
{\left(1\!-\!f_1\left(\varphi\left(e^{-ry+r\{\mu t+y\}}e^{-\nu s};0\right)\right)\right)\,dG_1(s)}\\
\!\!&-&\!\!(1-\alpha_1)e^{ry-r\{\mu t+y\}}\int\nolimits_t^{\infty}\!\!{\left(1\!-\!\varphi\!\left(e^{-ry+r\{\mu t+y\}}e^{-\nu s};0\right)\right)\,d\left(G_1\ast\overline{F}_{0,0}(s)\right)}.
\end{eqnarray*}

It follows from relations (\ref{I(t;u)_definition}) and (\ref{P(S(t)>x)_asymptotics}) that, for each $y\in\mathbb{R}$, there exists a number $\varepsilon_4\in(0,\mu-q{\sf E}Y^1)$ such that $e^{ry-r\{\mu t+y\}}I\left(t;\mu t+y\right)\leq C_{12}(y) e^{-t\Lambda(\mu-\varepsilon_4)}$, for some positive function $C_{12}(y)$, $y\in\mathbb{R}$, where $\Lambda(\mu-\varepsilon_4)>0$. Consequently, on account of the mean value theorem, applied to $f_1$, and the boundedness of the function $\left(1-\varphi(\lambda;0)\right)/\lambda$, $\lambda\geq0$, by some constant $C_{13}\geq c_{\ast}$, we have
\begin{eqnarray}
\!\!\!\!\!& &\!\!\!\!\!\left|I_2(t;y)\right|\leq C_{12}(y)e^{-t\Lambda(\mu-\varepsilon_4)}+\alpha_1 m_1\int\nolimits_t^{\infty}{\frac{1-\varphi(e^{-ry+r\{\mu t+y\}}e^{-\nu s};0)}{e^{-ry+r\{\mu t+y\}}e^{-\nu s}}e^{-\nu s}}\,d G_1(s)\label{I3leq_fluctuations}\\
\!\!\!\!\!&+&\!\!\!\!\!(1\!-\!\alpha_1)\!\!\!\int\limits_t^{\infty}{\!\frac{1\!-\!\varphi(e^{-ry+r\{\mu t+y\}}e^{-\nu s};0)}{e^{-ry+r\{\mu t+y\}}e^{-\nu s}}e^{-\nu s}\,d\left(G_1\ast\overline{F}_{0,0}(s)\right)}\!\leq\! C_{12}(y)e^{-t\Lambda(\mu-\varepsilon_4)}+C_{13}\!\left(1\!-\!\widetilde{G}(t)\right)\!\nonumber.
\end{eqnarray}
Here $\widetilde{G}$ is a c.d.f. such that $d\widetilde{G}(s)=e^{-\nu s}\,dG(s)$, $s\geq0$.

Let $t>2T$. Then by virtue of (\ref{limlimK(t,lambda)=0_fluctuations}), the mean value theorem applied to function $f_1$ and Lemma~\ref{L:lower_estimate_fluctuations} we obtain (for some positive constant $C_{14}$) the following relation
\begin{eqnarray}\label{I22leq_fluctuations}
\left|I_{12}(t,T;y)\right|&\leq&e^{ry-r\{\mu t+y\}}\int\nolimits_{t-T}^t{\left|E\left(t-s;\mu t+y\right)-1+\varphi(e^{-ry+r\{\mu t+y\}}e^{-\nu s};0)\right|\,d G(s)}\nonumber\\
&=&\int\nolimits_{t-T}^t{\left|K(t-s;y+\mu s)\right|e^{-\nu s}\,d G(s)}\leq C_{14}\left(1-\widetilde{G}(T)\right).
\end{eqnarray}

Applying the mean value theorem for $f_1$ once again we conclude that, for any $t>T$,
\begin{eqnarray}\label{I12leq_fluctuations}
& &\left|I_{11}(t,T;y)\right|\leq e^{ry-r\{\mu t+y\}}\int\nolimits_0^{t-T}{\left|E\left(t-s;\mu t+y\right)-1+\varphi(e^{-ry+r\{\mu t+y\}}e^{-\nu s};0)\right|\,d G(s)}\nonumber\\
&=&\int\nolimits_0^{t-T}{\left|K(t-s;y+\mu s)\right|e^{-\nu s}\,d G(s)}\leq\int_0^{t-T}{K_T(y+\mu s)\,d\widetilde{G}(s)}\leq{\sf E}{K_T\left(y+\mu \zeta\right)},
\end{eqnarray}
where $\zeta$ is a random variable with c.d.f. $\widetilde{G}$.

Combination of formulae (\ref{K(t;lambda)leqI2+I3_fluctuations})--(\ref{I12leq_fluctuations}) when $t>2T$ leads to the inequality
$$
\left|K(t;y)\right|\leq C_{12}(y)e^{-t\Lambda(\mu-\varepsilon_4)}+C_{13}\!\left(1\!-\!\widetilde{G}(t)\right)+C_{14}\left(1-\widetilde{G}(T)\right)+{\sf E}{K_T\left(y+\mu \zeta\right)}.
$$
It means that
$$
K_{2T}(y)\leq{\sf E}{K_T\left(y+\mu \zeta\right)}+C_{12}(y)e^{-T\Lambda(\mu-\varepsilon_4)}+\left(C_{13}+C_{14}\right)\left(1-\widetilde{G}(T)\right).
$$
According to the Lebesgue bounded convergence theorem, when $T\to\infty$ the latter formula entails the relation
$
K(y)\leq{\sf E}{K\left(y+\mu \zeta\right)}.
$
Iteration of this formula yields the following inequality
\begin{equation}\label{K(lambda)leq_fluctuations}
K(y)\leq{\sf E}{K\left(y+\mu Z_n\right)},
\end{equation}
where $Z_n:=\sum_{k=1}^{n}\zeta_k$ and $\zeta_k$, $k\in\mathbb{Z}^+$, are independent identically distributed random variables with the same distribution as $\zeta$. On account of the strong law of large numbers and the Lebesgue bounded convergence theorem inequality (\ref{K(lambda)leq_fluctuations}) implies
$
0\leq K(y)\leq K(+\infty).
$
However, $K(+\infty)=0$ in view of (\ref{limlimK(t,lambda)=0_fluctuations}). Thus, Lemma~\ref{L:limK=0_fluctuations} is proved. $\square$

\vskip0.2cm
\noindent\emph{Proof of Corollary~\emph{\ref{C:d=1}}.} For $N=1$ and $x=w_1=0$, the statement of Corollary~\ref{C:d=1} is assured by Lemma~\ref{L:limK=0_fluctuations}. Now we deal with $N>1$ and $x\in W$, say, $x=w_i$. Discuss here the main differences in the cases of a single catalyst and several catalysts and give the subsequent proof omitting cumbersome details. In the setting of the problem with several catalysts a counterpart of equation (\ref{E(t;u)_equation}) in Lemma~\ref{L:equation_multi} is a system of integral equations
\begin{eqnarray}
& &{\sf P}_{w_i}\left(M_t>u\right)=\alpha_i\int\nolimits_0^t{\left(1-f_i\left(1-{\sf P}_{w_i}\left(M_{t-s}>u\right)\right)\right)\,dG_i(s)}\label{NP_0(Ltgamma/Mtgamma<lambda)_equation_fluctuations}\\
&+&(1-\alpha_i)\sum_{j=1}^N\int\nolimits_0^t{{\sf P}_{w_j}\left(M_{t-s}>u\right)\,d\left(G_i\ast {_{W_j}\overline{F}_{w_i,w_j}(s)}\right)}+I_i^{(N)}(t;u),\nonumber
\end{eqnarray}
where $i=1,\ldots,N$ and functions $I_i^{(N)}(t;u)$, $t\geq0$, $u\geq\max\{w_1,\ldots,w_N\}$, are of the form
\begin{eqnarray}\label{NJ-1(t;a)=_fluctuations}
\!\!\!\!\!\!\!\!\!& &\!\!I_i^{(N)}(t;u)={\sf P}_{w_i}\left(S(t)>u\right)-\sum_{k=1}^N\int\nolimits_0^t{{\sf P}_{w_k}\left(S(t-s)>u\right)\,d\,{_{W_k}F_{w_i,w_k}(s)}}\nonumber\\
\!\!\!\!\!\!\!\!\!&-&\!\!\alpha_i\int\nolimits_0^t{{\sf P}_{w_i}\left(S(t-s)>u\right)\,d G_i(s)}+\alpha_i\sum_{k=1}^N\int\nolimits_0^t{{\sf P}_{w_k}\left(S(t-s)>u\right)\,d G_i\ast{_{W_k}F_{w_i,w_k}(s)}}.
\end{eqnarray}
The next step in the case of several catalysts is to introduce a counterpart of function $G$ appearing in Lemma~\ref{L:fluctuations_main_step}, namely, the matrix $\mathcal{G}(t)=\left(G^{(N)}_{i,j}(t)\right)_{i,j=1}^N$, where $G^{(N)}_{i,j}(t):=\delta_{i,j}\alpha_i m_i G_i(t)+(1-\alpha_i)G_i\ast{_{W_j}\overline{F}_{w_i,w_j}(t)}$, $t\geq0$, and $\delta_{i,j}$ is the Kronecker delta. Note that an entry $d_{i,j}(\lambda)$ of the matrix $D(\lambda)$, $\lambda\geq0$, arising in the definition of the supercritical regime of CBRW (see, e.g., \cite{B_TVP_14}), is the Laplace transform of function $G^{(N)}_{i,j}$.

Let us turn to a counterpart of Lemma~\ref{L:lower_estimate_fluctuations} and afterwards return to counterparts of Lemmas~\ref{L:fluctuations_main_step} and \ref{L:int<=}. In accordance with the mean value theorem applied to functions $f_1,\ldots,f_N$, the equations system (\ref{NP_0(Ltgamma/Mtgamma<lambda)_equation_fluctuations}) entails the following vector inequality, valid coordinate-wise,
\begin{equation}\label{NP_0(Ltgamma/Mtgamma<lambda)_inequality_fluctuations}
\mathcal{P}(t;u)\leq\mathcal{G}\ast\mathcal{P}(t;u)+\mathcal{I}(t;u),
\end{equation}
where $\mathcal{P}(t;u):=\left({\sf P}_{w_1}\left(M_t>u\right),\ldots,{\sf P}_{w_N}\left(M_t>u\right)\right)^{\top}$ and $\mathcal{I}(t;u):=\left(I^{(N)}_1(t;u),\ldots,I^{(N)}_N(t;u)\right)^{\top}$, and the symbol $\top$ denotes the matrix transposition. Recall that operation ``$\ast$'' of the matrix convolution is defined in the same manner as the matrix multiplication except for that we convolve the entries rather than multiply them. Iterating inequality (\ref{NP_0(Ltgamma/Mtgamma<lambda)_inequality_fluctuations}) $k$ times, letting $k$ tend to infinity and employing Lemma~1.1 in paper \cite{Crump_70}, similarly to the proof of formula (\ref{P_0(Ltgamma/Mtgamma<lambda)_inequality_fluctuations}), we deduce that
$$
\mathcal{P}(t;u)\leq\sum_{k=0}^{\infty}\mathcal{G}^{\ast k}\ast\mathcal{I}(t;u).
$$
Thus, alike to Lemma~\ref{L:fluctuations_main_step} for $N=1$, in the case of $N>1$ we investigate the asymptotic behavior of the expression $\sum_{k=0}^{\infty}\mathcal{G}^{\ast k}\ast\mathcal{I}(t;u)$, when $u=\mu t+y$ and $t\to\infty$. Completely similarly to Lemma~\ref{L:fluctuations_main_step}, using Corollary~3.1, item (i), of paper \cite{Crump_70} (instead of Theorem~25 in \cite{Vatutin_book_09}, p.~30, and results on p.~55 in \cite{Cox_62}), we conclude that
$$
\sum_{k=0}^{\infty}\mathcal{G}^{\ast k}\ast\mathcal{I}(t;\mu t+y)\sim e^{-ry+r\{\mu t+y\}}\left(K_1^{(N)},\ldots,K^{(N)}_N\right)^{\top},\quad t\to\infty.
$$
The constants $K_i^{(N)}>0$, $i=1,\ldots,N$, can be written in an explicit form which is cumbersome and superfluous and, therefore, is omitted. Moreover, Lemmas~\ref{L:int<=} and \ref{L:lower_estimate_fluctuations} remain unchanged as well in the case $N>1$ (with, possibly, other constants $C\,'_4$ and $C\,'$ instead of $C_4$ and $C$, respectively).

The generalization of the function $J(t;u)$, $t\geq0$, $u\in\mathbb{R}$, in the case of $N>1$ is a vector-valued function $\mathcal{J}(t;u)$, $t\geq0$, $u\in\mathbb{R}$, with coordinates $J^{(N)}_i(t;u)$, $i=1,\ldots,N$, of the form
$$
m_i\int\nolimits_0^t{{\sf P}_{w_i}\left(M_{t-s}>u\right)\,dG_i(s)}-\int\nolimits_0^t{\left(1-f_i\left(1-{\sf P}_{w_i}\left(M_{t-s}>u\right)\right)\right)\,dG_i(s)}.
$$
A counterpart of Lemma~\ref{L:I-1(t;lambdaLt)ast_sim_fluctuations} in the case of several catalysts asserts that under the same conditions one has
$$
\lim_{y\to+\infty}\lim_{t\to\infty}e^{ry}\sum_{k=0}^{\infty}\mathcal{G}^{\ast k}\ast\mathcal{J}(t;\mu t+y)=(0,\ldots,0)^{\top}.
$$
The proof repeats the proof of Lemma~\ref{L:I-1(t;lambdaLt)ast_sim_fluctuations}, although now we apply Corollary~3.1, item (i), of paper \cite{Crump_70} instead of Theorem~25 in book \cite{Vatutin_book_09}, p.~30.

The differences between the statements of Lemmas~\ref{L:limlim=theta_fluctuations}, \ref{L:limK=0_fluctuations} and their corresponding counterparts in the case of several catalysts are insignificant, the same observation refers to their proofs. Hence, we only note that the proof of a counterpart of Lemma~\ref{L:limK=0_fluctuations} follows the proof of Theorem~3.3 of paper \cite{Kaplan_75}, whereas the proof of Lemma~\ref{L:limK=0_fluctuations} is based on the work~\cite{Athreya_69}. Thus, Corollary~\ref{C:d=1} is established in the case of $N\geq1$ and the starting point $x\in W$.

It remains to prove Corollary~\ref{C:d=1} in the case of $N\geq1$ and $x\notin W$. The case of the starting point $x\notin W$ is reduced to the case of $N+1$ catalysts, since we may set $w_{N+1}=x$, $\alpha_{N+1}=0$, $m_{N+1}=0$,
$G_{N+1}(t)=1-e^{-qt}$, $t\geq0$. According to Lemma~3 in paper \cite{B_TVP_14} the new CBRW with the catalysts set $\{w_1,\ldots,w_{N+1}\}$ is supercritical, whenever the initial CBRW is supercritical, and the Malthusian parameters in these models coincide. Therefore, we may employ the proven part of Corollary~\ref{C:d=1} for the case of $N+1$ catalysts and the starting point from the set $W$. Corollary~\ref{C:d=1} is proved completely. $\square$

\vskip0.2cm
\noindent\emph{Proof of Theorem~\emph{\ref{T:fluctuation_light}}.} When $d=1$ and $N\in\mathbb{N}$ the statement of Theorem~\ref{T:fluctuation_light} coincides with that of Corollary~\ref{C:d=1}. Consider the case $d>1$ and $N\in\mathbb{N}$. The main traits of the proof of Theorem~\ref{T:fluctuation_light} are the same as those of Corollary~\ref{C:d=1}. Therefore, discuss the main differences only. Firstly, a counterpart of the integral equation in Lemma~\ref{L:equation_multi} in the case of CBRW on $\mathbb{Z}^d$ with $N$ catalysts has the following form
\begin{eqnarray*}
& &E^{\,{\bf r}}_{{\bf w}_i}(t;u)=\alpha_i\int\nolimits_0^t{\left(1-f_i\left(1-E^{\,{\bf r}}_{{\bf w}_i}(t-s;u)\right)\right)\,dG_i(s)}\\
&+&(1-\alpha_i)\sum_{j=1}^N\int\nolimits_0^t{E^{\,{\bf r}}_{{\bf w}_j}(t-s;u)\,d\left(G_i\ast {_{W_j}\overline{F}_{w_i,w_j}(s)}\right)}+I_i^{\,{\bf r},N}(t;u),
\end{eqnarray*}
where $E^{\,{\bf r}}_{{\bf w}_i}(t;u):={\sf P}_{{\bf w}_i}\left(\exists z\in Z(t):\langle{\bf X}^z(t),{\bf r}\rangle>u\right)={\sf P}_{{\bf w}_i}\left(M_t({\bf r})>u\right)$, $t\geq0$, $u\geq\max\{\langle{\bf w}_j,{\bf r}\rangle:j=1,\ldots,N\}$, $i=1,\ldots,N$, and the function $I_i^{\,{\bf r},N}$ coincides with the function $I_i^{(N)}$ except to the replacement of expression $S(t)$ by $\langle{\bf S}(t),{\bf r}\rangle$ in its definition.

Secondly, in a counterpart of Lemma~\ref{L:P(S(t)>x)_from_discrete_to_continuous} we now consider a random walk $\left\{\langle{\bf S}(t),{\bf r}\rangle,t\geq0\right\}$ instead of the random walk $\left\{S(t),t\geq0\right\}$. As noted in Section~\ref{s:results_fluctuations}, the random variables $S(t)$ have lattice (arithmetical) distribution, for each $t\geq0$, whereas all the random variables in the set $\langle{\bf S}(t),{\bf r}\rangle$, $t\geq0$, may have either a lattice distribution or a non-lattice one. Hence, for non-lattice distributions, we use formula (6.1.16) from Corollary~6.1.7 in book \cite{Borovkov_Borovkov_08}, Ch.~6, Section~1, instead of formula (6.1.17) from the same corollary. Namely, a counterpart of relation (\ref{P(S(t)>x)_asymptotics}) in Lemma~\ref{L:P(S(t)>x)_from_discrete_to_continuous} in the case of non-lattice distribution for $d>1$ and $N\in\mathbb{N}$ is the following relation
$$
{\sf P}_{\bf 0}\left(\langle{\bf S}(t),{\bf r}\rangle\geq x\right)\sim\frac{e^{-t\Lambda_{\bf r}(\theta)}}{\lambda_{\bf r}(\theta)\sqrt{2\pi t D_{\bf r}(\theta)}},
$$
as $x,t\to\infty$, uniformly in $\theta=x/t\in[q\langle{\sf E}{\bf Y}^1,{\bf r}\rangle+\varepsilon_1,\Theta_1]$, $x\in\mathbb{R}$, $t>0$. Here $\Lambda_{\bf r}(\vartheta):=\sup_{s\in\mathbb{R}}\left(\vartheta s-\ln{\sf E}e^{s\langle{\bf S}(1),{\bf r}\rangle}\right)=\sup_{s\in\mathbb{R}}\left(\vartheta s-H_{\bf r}(s)\right)$, $\lambda_{\bf r}(\vartheta):=\Lambda_{\bf r}'(\vartheta)$, $D_{\bf r}(\vartheta):=\left.H_{\bf r}''(s)\right|_{s=\lambda_{\bf r}(\vartheta)}$, $\vartheta\in\mathbb{R}$, and $H_{\bf r}(s):=H(s{\bf r})$, $s\in\mathbb{R}$. The subsequent argument proving Corollary~1 easily extends to the general case, provided that additionally one replaces $\mu$ by $\nu$ and $r$ by $1$. Obviously, now ${\sf E}e^{s\langle{\bf S}(t),{\bf r}\rangle}=e^{tH_{\bf r}(s)}$, $\theta_0=H_{\bf r}'(1)$,
$\theta_0>H_{\bf r}(1)=\nu$. The further details can be omitted.

Thus, Theorem~\ref{T:fluctuation_light} is proved completely. $\square$

\vskip0.2cm The author expresses acknowledgements to Professor V.A.~Vatutin and Professor S.G.~Foss for advice and useful discussions.

\end{document}